\documentclass[11pt]{article}
\usepackage{amsmath}
\usepackage{mathrsfs}
       \usepackage{amsfonts}
       \usepackage{stmaryrd}
       \usepackage{latexsym,amssymb,mathrsfs,fancyhdr}
       \usepackage[square, comma, sort&compress, numbers]{natbib}
       \font\tenmsb=msbm10
       \font\sevenmsb=msbm7
       \font\fivemsb=msbm5
       \catcode`\@=11
       \ifx\amstexloaded@\relax\catcode`\@=\active
       \endinput\else\let\amstexloaded@\relax\fi
       \def\spaces@{\space\space\space\space\space}
       \def\spaces@@{\spaces@\spaces@\spaces@\spaces@\spaces@}
       \def\space@.  {\futurelet\space@\relax}
       \space@.   %
       \def\Err@#1{\errhelp\defaulthelp@\errmessage{AmS-TeX error: #1}}
       \def\relaxnext@{\let\next\relax}
       \def\accentfam@{7}
       \def\noaccents@{\def\accentfam@{0}}
       \def\Cal{\relaxnext@\ifmmode\let\next\Cal@\else
       \def\next{\Err@{Use \string\Cal\space only in math mode}}\fi\next}
       \def\Cal@#1{{\Cal@@{#1}}}
       \def\Cal@@#1{\noaccents@\fam\tw@#1}
       \def\Bbb{\relaxnext@\ifmmode\let\next\Bbb@\else
       \def\next{\Err@{Use \string\Bbb\space only in math mode}}\fi\next}
       \def\Bbb@#1{{\Bbb@@{#1}}}
       \def\Bbb@@#1{\noaccents@\fam\msbfam#1}
       \newfam\msbfam
       \textfont\msbfam=\tenmsb
       \scriptfont\msbfam=\sevenmsb
       \scriptscriptfont\msbfam=\fivemsb

\newtheorem{Remark}{Remark}[section]

\@addtoreset{equation}{section}

\newcommand{\beq}{\begin{equation} }
\newcommand{\eeq}{\end{equation} }

\usepackage[english]{babel}
\usepackage{amsmath,amssymb}

\newtheorem{lemma}{Lemma}[section]
\newtheorem{iteration lemma}{iteration Lemma}[section]

\begin{document}

\setlength{\columnsep}{5pt}
\title{Multiple semiclassical states for  fractional Schr\"{o}dinger equations with asymptotically linear nonlinearities\footnote{The work was supported by the National Natural Science Foundation of China (No.
11871146) and Qinglan Project of Jiangsu Province of China.}}

\author{{ Hui Zhang$^{1,2}$\footnote{Corresponding author. Email addresses: huihz0517@126.com (H. Zhang), zhangfubao@seu.edu.cn (F.B. Zhang).}, Fubao Zhang$^{3}$}\\1.~ Department of  Mathematics, Jinling Institute of Technology, \\ Nanjing 211169, China
\\2.~ Department of  Mathematics,
 Nanjing University, Nanjing 210093, China\\
 3.~ Department of  Mathematics,
 Southeast University, Nanjing 210096, China}
\date{}
\maketitle
{\bf Abstract}: {In this paper, we consider the singularly perturbed
fractional Schr\"{o}dinger  equation
\begin{equation*}
\epsilon^{2\alpha}(-\Delta)^\alpha u+V(x)u=f(u),\quad x\in \mathbb{R}^N,
\end{equation*}
where $\epsilon>0$ is a small parameter, $\alpha\in(0,1)$, $(-\Delta)^\alpha$ is  the fractional Laplacian operator of order $\alpha$, $V$ possesses global minimum points, and $f$ is asymptotically linear at infinity. We investigate the relationship between the number of positive solutions and the topology of the set where the potential $V$ attains its global minimum. We also construct multiple concentrating solutions if $V$ has several strict global minimum points. In particular, some new tricks and the method of Nehari manifold dependent on a suitable restricted set are introduced to overcome the difficulty resulting from the appearance of asymptotically linear nonlinearity.}



\section{Introduction and main results}
\setcounter{equation}{0}

In the last few years, a great deal of work has been devoted to the study of fractional  Schr\"{o}dinger equation
\begin{equation}\label{1.1.1}i\epsilon\frac{\partial\Psi}{\partial t}=\epsilon^{2\alpha}(-\Delta)^\alpha\Psi+W(x)\Psi-g(|\Psi|)\Psi,\quad (t,x)\in\mathbb{R}\times\mathbb{R}^N, \end{equation}
where $\epsilon$ is a positive constant, $i$ is the imaginary unit, $\alpha\in(0,1)$, $(-\Delta)^\alpha$ is  the fractional Laplacian operator of order $\alpha$, $W$ is the potential, and $g$ is a nonlinear term.
 Equation (\ref{1.1.1}) was introduced by Laskin \cite{LN1,LN2}, and comes from an expansion of the Feynman path integral from Brownian-like
to L\'{e}vy-like quantum mechanical paths.
 It also appeared
in several areas such as optimization, finance, phase transitions, stratified materials, crystal dislocation, flame propagation, conservation laws. We refer to \cite{CT,WZ} for more background.
The standing wave solution is a  solution of the form $\Psi(x,t)=u(x)e^{\frac{-iEt}{\hbar}}$. Then $\Psi(x,t)$ is a solution of (\ref{1.1.1}) if and only if $u$ is a solution of the equation
\begin{equation}\label{1.1}\epsilon^{2\alpha}(-\Delta)^\alpha u+V(x)u=f(u), \quad \ x\in\mathbb{R}^N,\end{equation}
where $V(x)=W(x)-E$ and $f(u)=g(|u|)u$. In this paper we are interested in nontrivial solutions of (\ref{1.1}) for small $\epsilon>0$. For small $\epsilon>0$, standing wave solutions are called semi-classical states. The limit
$\epsilon\rightarrow0$ describes the transition from quantum mechanics to classical mechanics.

In the case that $\alpha=1$, (\ref{1.1}) turns out  to be the Schr\"{o}dinger equation
 \begin{equation}\label{1.2.0}-\epsilon^{2}\Delta u+V(x)u=f(u), \quad \ x\in\mathbb{R}^N,\end{equation}and there are large amounts of works on the  existence, multiplicity and concentration of solutions for (\ref{1.2.0}) under various conditions on $V$ and $f$. In a well-known paper \cite{Rap}, Rabinowitz firstly introduced the global condition\vskip 0.1 true cm
\noindent (V)\ \  $\liminf_{|x|\rightarrow\infty}V(x)>\inf_{\mathbb{R}^N}V(x)>0,$\\
 and used Mountain Pass theorem to show the equation
\begin{equation}\label{1.2}-\epsilon^2\Delta u+V(x)u=|u|^{p-2}u,\quad x\in\mathbb{R}^N,\end{equation}
with $2<p<2^*=2N/(N-2)$ if $N\geq3$, $2<p<\infty$ if $N=1,2$,
 possesses a
positive ground state for small $\epsilon>0$.
Moreover, Wang \cite{WXF} proved that positive ground states of (\ref{1.2}) concentrate around global minimum points of $V$ as $\epsilon\rightarrow0$. Later, under the same condition (V), Cingolani and  Lazzo \cite{CL} investigated the relationship between the number of positive solutions
of (\ref{1.2}) and the topology of the set where $V$ attains its global minimum using Ljusternik-Schnirelmann category theory.
Cao and Noussair \cite{CN}  studied the
problem
\begin{equation}\label{1.1.4}-\Delta u+\mu u=Q(x)|u|^{p-2}u, \quad  x\in \mathbb{R}^N,\end{equation}
where $2<p<2^*$, $N\geq3$ and $\mu>0$. By assuming $Q\in C(\mathbb{R}^N)$ and
   \vskip 0.1 true cm
\noindent (Q)\ $Q(x)\geq0$ in $\mathbb{R}^N$ and there exist $a^1,...,a^k$ in $\mathbb{R}^N$ such that $Q(a^j)$ are strict maximums and satisfy $$Q(a^j)=Q_{max}:=\max_{x\in \mathbb{R}^N}Q(x)>0, \ j=1,...,k,$$
Cao and Noussair proved that there exists $\mu_0>0$ such that (\ref{1.1.4}) has at least  $k$ positive and nodal solutions when $\mu\geq\mu_0$. Their results showed how the properties of $Q(x)$ affects the number of solutions. Subsequently, the idea of \cite{CN} was adapted by Chen et al. \cite{CKW} to deal with
singularly perturbed 4-superlinear  Schr\"{o}dinger-Poisson system.


In the case that $\alpha\in (0,1)$, the equation (\ref{1.1}) was widely studied and many researchers focused on superlinear nonlinearity.  Fall et al. \cite{FMV}
 considered the problem
\begin{equation}\label{1.2.3}\epsilon^{2\alpha}(-\Delta)^\alpha u+V(x)u=|u|^{p-2}u, \ \ \ \ x\in \mathbb{R}^N,\end{equation}where $2<p<2^*_\alpha:=\frac{2N}{N-2\alpha}$, and showed that if the potential $V$ is coercive and has a unique global minimum point, then the ground state of (\ref{1.2.3}) concentrates at such a minimum point as $\epsilon\rightarrow0$. When $V$ is sufficiently smooth, D\'{a}vila et al. in \cite{DPW} recovered various existence results known
for the case $\alpha=1$ and showed the existence of solutions around $k$ nondegenerate
critical points of $V$ for (\ref{1.2.3}). In \cite{DMS,SZ4}, the authors dealt with (\ref{1.1}) with vanishing potentials, and established the existence, multiplicity and concentration of solutions using variational methods.
By means of penalization method, Alves and Miyagaki \cite{AM} proved the existence and
concentration of positive solutions of (\ref{1.1})
 when $V$ possesses local minimum points. When $f$ is superlinear and subcritical,  Figueiredo and Siciliano \cite{GF} showed the existence and multiplicity of positive solutions of (\ref{1.1}) by using Ljusternik-Schnirelmann category theory and
Morse theory, and exploiting the topological complexity of the set of global minimum points of the potential $V$. Regarding the critical growth case, Shang and Zhang \cite{SZ2} considered  the equation
\begin{equation}\label{1.2.4}
\epsilon^{2\alpha}(-\Delta)^\alpha u+V(x)u=|u|^{2^*_\alpha-2}u+\lambda f(u),\ x\in \mathbb{R}^N,\end{equation}
where $f$ is a subcritical and superlinear term, proved that (\ref{1.2.4})
has a ground state and investigated the relation between the number of solutions and the topology of the set where $V$ attains its global minimum for large $\lambda$ and small $\epsilon$. Later, for the equation (\ref{1.2.4}) with $\lambda=1$, under a local condition imposed on
the potential $V$, He and Zou \cite{HZ1} showed
the multiplicity of positive solutions concentrating around local minimum points
of $V$. Shang and Zhang \cite{SZ3} extended the results in \cite{CN} about the equation (\ref{1.1.4}) to (\ref{1.2.4}) with superlinear nonlinearity $f$.

On the contrary, the results about semiclassical states of (\ref{1.1}) with asymptotically linear $f$ are relatively few. When the nonlinearity $f$ is either asymptotically linear or superlinear, Ambrosio \cite{Ambro} proved  that (\ref{1.1}) has a positive solution which concentrates around the local minimum point of the potential $V$. In \cite{LLC}, Li et al. considered the equation (\ref{1.1}) with asymptotically linear nonlinearity and showed that (\ref{1.1}) possesses a ground state which concentrates near global minimum point of linear potential or  global maximum point of nonlinear potential.

To the best of our knowledge, when $f$ is asymptotically linear,  the multiplicity of concentrating solutions of the equation (\ref{1.1}) or (\ref{1.2.0}) has not been studied before. Consequently, the purpose of the paper is to fill this gap. Since the treatments of (\ref{1.1}) and (\ref{1.2.0}) are similar, in this paper we just consider the equation (\ref{1.1}) with $N\geq3$, and
the assumptions on $V$ and $f$ are as follows.
  \vskip 0.1 true cm
 \noindent(V$_1$) $V\in C(\mathbb{R}^N)\cap L^\infty(\mathbb{R}^N)$ and  $V_\infty:=\liminf_{|x|\rightarrow+\infty}V(x)>V_0:=\inf_{\mathbb{R}^N}V(x)>0$.
   \vskip 0.1 true cm
 \noindent(f$_1$) $f\in C^1(\mathbb{R})$, $f(t)=0$ for $t\leq0$, and $f(t)=o(t)$ as $t\rightarrow0$.
  \vskip 0.1 true cm
 \noindent(f$_2$) There exist $2<q<2^*_\alpha:=\frac{2N}{N-2\alpha}$ and $C_0>0$ such that $|f'(t)|\leq C_0(1+|t|^{q-2})$,\ $\forall t\in\mathbb{R}$.
 \vskip 0.1 true cm
 \noindent(f$_3$) $|V|_{L^\infty(\mathbb{R}^N)}<\lim_{t\rightarrow+\infty}\frac{f(t)}{t}=l_0<+\infty$.
  \vskip 0.1 true cm
 \noindent(f$_4$) $t\mapsto\frac{f(t)}{t}$ is increasing in $(0,+\infty)$.
  \vskip 0.1 true cm
  \noindent(f$_5$) $\lim_{t\rightarrow+\infty}\bigl[\frac12 f(t)t-F(t)\bigr]=+\infty$, where $F(t):=\int^t_0f(s)ds$.
    \vskip 0.1 true cm

\begin{Remark} A special case  of $f$ is the saturable nonlinearity $$f(u)=\frac{u^3}{1+su^2} \ \text{if}\ u>0, \ \text{and}\ f(u)=0\ \text{if}\ u\leq0,$$ where $0<s<|V|^{-1}_{L^\infty(\mathbb{R}^N)}$. The Schr\"{o}dinger equation with saturable nonlinearity arises in nonlinear optics (see \cite{KA}) and is our original motivation for studying (\ref{1.1}) with
asymptotically linear nonlinearity. \end{Remark}

To state the first result about the multiplicity and concentration of solutions, set
$$\mathcal{V}:=\{x\in\mathbb{R}^N:V(x)=V_0\},\ \mathcal{V}_\delta:=\{x\in\mathbb{R}^N:dist(x,\mathcal{V})\leq\delta\},$$
for $\delta>0$. Without loss of generality, below we assume $0\in \mathcal{V}$. Recall that, if
$Y$ is a closed subset of a topological space $X$, the Ljusternik-Schnirelmann category $cat_X(Y)$ is the least number of closed and contractible sets in $X$ which cover $Y$.

\textbf{Theorem 1.1.} Let (V$_1$) and (f$_1$)-(f$_5$) hold. Then
 for each $\delta>0$, there exists $\epsilon_\delta>0$ such that for any $\epsilon\in(0,\epsilon_\delta)$, the equation (\ref{1.1}) has at least $cat_{\mathcal{V}_\delta}(\mathcal{V})$ positive solutions in $H^\alpha(\mathbb{R}^N)$. Moreover, if $u_\epsilon$ denotes one of these positive solutions and $\eta_\epsilon\in\mathbb{R}^N$ its global maximum point, then
 $$\lim_{\epsilon\rightarrow0}V(\eta_\epsilon)=V_0,$$
 and
 $$u_\epsilon(x)\leq\frac{C\epsilon^{N+2\alpha}}{\epsilon^{N+2\alpha}+|x-\eta_\epsilon|^{N+2\alpha}},$$
for some constant $C>0$.
  \vskip 0.2 true cm
Next we state another result on the multiplicity of concentrating solutions and assume

\vskip 0.1 true cm
\noindent(V$_2$) $a^1,...,a^k\in \mathcal{V}$, and each of $a^1,...,a^k$ is a strict global minimum point of $V$.
  \vskip 0.2 true cm
\textbf{Theorem 1.2.} Assume that (V$_1$), (V$_2$) and (f$_1$)-(f$_5$) are satisfied. Then there exists $\epsilon^*>0$ such that for every $\epsilon\in (0,\epsilon^*)$,

\noindent(1) problem (\ref{1.1}) admits at least $k$ positive solutions $w^j_\epsilon$ in $H^\alpha(\mathbb{R}^N)$, $j=1,...,k$;

\noindent(2) if additionally $V$ is uniformly continuous, then each $w^j_\epsilon$ satisfies

\noindent(i) there exists a maximum point $x^j_\epsilon\in\mathbb{R}^N$ of $w^j_\epsilon$ such that $x^j_\epsilon\rightarrow a^j$ as $\epsilon\rightarrow0$. Setting $v^j_\epsilon(x):=w^j_\epsilon(\epsilon x+x^j_\epsilon)$, then $v^j_\epsilon$ converges strongly in $H^\alpha(\mathbb{R}^N)$ to a positive solution $v^j$ of the equation
\begin{equation}\label{1.0.0}
 (-\Delta)^\alpha u+V_0u=f(u), \quad \text{in}\ \mathbb{R}^N;
\end{equation}

\noindent(ii)$$w^j_\epsilon(x)\leq\frac{C_j\epsilon^{N+2\alpha}}{\epsilon^{N+2\alpha}+|x-x^j_\epsilon|^{N+2\alpha}},$$
for some constants $C_j>0$.

  \vskip 0.1 true cm
\begin{Remark}In Theorem 1.2, we know that the
global maximum point of each positive solution concentrates around the associated global minimum
point of the potential $V$, which makes concentration phenomenon become more specific than in Theorem 1.1.
\end{Remark}

Now we state the third result which implies that the problem (\ref{1.1}) has at least one ground state among the solutions of Theorem 1.2 if the potential $V$ has no other global minimum points than $a^1$,...,$a^k$. Suppose that

  \vskip 0.1 true cm
\noindent(V$_3$) $\mathcal{V}=\{a^1,...,a^k\}$, where $a^1,...,a^k$ are given in (V$_2$).

    \vskip 0.2 true cm
\textbf{Theorem 1.3.} Let (V$_1$)-(V$_3$) and (f$_1$)-(f$_5$) hold. Then there exists $\bar{{\epsilon}}>0$ such that for each
$\epsilon\in (0,\bar{{\epsilon}})$, we can find at least a ground state of (\ref{1.1}) among the solutions of Theorem 1.2.

\begin{Remark} It is worth pointing out that there is no related result for (\ref{1.2.0}) with asymptotically linear nonlinearity. So the results are also new for (\ref{1.2.0}). Namely, under the conditions (V$_1$)-(V$_3$) and (f$_1$)-(f$_5$), Theorems 1.1, 1.2 and 1.3 can be extended trivially to the Schr\"{o}dinger equation (\ref{1.2.0})
except for changing polynomial decay of solutions into exponential decay of solutions.
\end{Remark}

\begin{Remark} We shall use the method of Nehari manifold dependent on a restricted set, similar idea has been used to show the existence of a ground state or a sign-changing solution for some elliptic equations, see \cite{ZZJMAA,ZT} for example, but as we know, this is the first time to use this idea to seek multiple concentrating solutions. \end{Remark}

To a certain extent, our results can be viewed as the counterparts of the results for (\ref{1.1}) with superlinear nonlinearity in \cite{GF,SZ2,SZ3} to asymptotically linear case. Due to asymptotically linear nonlinearity, we cannot use variational methods in a standard way. The main difficulties lie in three aspects: (i) the asymptotically linear nonlinearity leads that
the method of  Nehari manifold can not be used directly.
(ii) The asymptotically linear nonlinearity $f$ makes it complex to show the
PS condition holds, which is essential for using Ljusternik-Schnirelmann category theory to find multiple solutions.
(iii) Since  $f$ is asymptotically linear, it is difficult to obtain the strong convergence of positive solutions as $\epsilon\rightarrow0$.
Concerning (i) and (ii), we find a suitable restricted set \begin{equation}\label{1.3}\Theta_\epsilon:=\{u\in H^\alpha(\mathbb{R}^N): [u]^2_\alpha+\int_{\mathbb{R}^N}V(\epsilon x)u^2dx-l_0|u|^2_2<0\},\end{equation}
where $[\cdot]_\alpha$ and $l_0$ are  given in (\ref{2.1.0}) and (f$_3$) respectively. Then the method of Nehari manifold dependent on $\Theta_\epsilon$ can be developed and we use it to show the multiplicity and concentration of solutions.
In this process, on one hand, the  least energy  has a good characterization dependent on  $\Theta_\epsilon$ which plays an important role in recovering the compactness of PS sequence. On the other hand,  the
involvement of the set $\Theta_\epsilon$ in the method of Nehari manifold enforces the implementation of some new estimates about $\Theta_\epsilon$, and then some verifications turn out to be difficult. Especially, in the proof of Theorem 1.2, it is complicated to verify that the subset of Nehari manifold associated with the strict global minimum point of the potential $V$ is non-empty. As for (iii), by virtue of a similar argument as the global compactness lemma, we deduce that positive solutions of (\ref{1.1}) converge strongly to those of the limit equation.

The paper is organized as
follows. In Section 2 we introduce the variational framework. In Section 3 we study the autonomous problem. In Sections 4, 5 and 6, we  give the proof of Theorems
1.1, 1.2 and 1.3 respectively.
\section{Preliminaries and functional setting}

In this paper, we use the following
notations.  For $1\leq p\leq\infty$, the
norm in $L^p(\mathbb{R}^N)$ is denoted by $|\cdot|_{p}.$ For any $r>0$ and $x\in\mathbb{R}^N$,
$B_r(x)$ denotes the ball  centered at $x$ with the radius $r.$
$H^\alpha(\mathbb{R}^N)$ is the fractional Sobolev space  defined by
$$H^\alpha(\mathbb{R}^N)=\bigl\{u\in L^2(\mathbb{R}^N):\frac{|u(x)-u(y)|}{|x-y|^{\frac{N}{2}+\alpha}}\in L^2(\mathbb{R}^N\times\mathbb{R}^N)\bigr\},$$
and endowed with the standard norm
$$\|u\|
=\Bigl([u]^2_\alpha+\int_{\mathbb{R}^N}|u|^2dx
\Bigr)^{\frac{1}{2}},
$$
while
\begin{equation}\label{2.1.0}[u]^2_\alpha:=\int_{\mathbb{R}^N}\int_{\mathbb{R}^N}\frac{|u(x)-u(y)|^2}{|x-y|^{N+2\alpha}}dxdy,\end{equation}
is the Gagliardo (semi) norm.

 By making the change of variable $x\rightarrow\epsilon x$, the problem (\ref{1.1}) turns out to be
\begin{equation}\label{2.1}
 (-\Delta)^\alpha u+V(\epsilon x)u=f(u),\quad x\in
\mathbb{R}^{N}.
\end{equation}
By (V$_1$) we have
$$\|u\|_\epsilon=\Bigl([u]^2_\alpha+\int_{\mathbb{R}^N}V(\epsilon x)u^2dx\Bigr)^{\frac12},$$
is an equivalent norm on $H^\alpha(\mathbb{R}^N)$.
The functional associated with the equation (\ref{2.1}) is$$
I_\epsilon(u)=\frac{1}{2}\|u\|^2_\epsilon-\int_{\mathbb{R}^N}F(u)dx.$$

For the readers' convenience, we recall some results of the space $H^\alpha(\mathbb{R}^N)$.
\begin{lemma}{\label{l2.0.1}}(\cite{DP}) The embedding $H^\alpha(\mathbb{R}^N)\hookrightarrow L^q(\mathbb{R}^N)$ is continuous for any $q\in [2,2^*_\alpha]$. Moreover, the embedding $H^\alpha(\mathbb{R}^N)\hookrightarrow L^q(\mathbb{R}^N)$ is locally compact whenever $q\in[1,2^*_\alpha)$.
\end{lemma}

\begin{lemma}(\cite[Lemma 2.4]{SS1}){\label{l2.0.2}} Assume $\{u_n\}$ is bounded in $H^\alpha(\mathbb{R}^N)$, and it satisfies
\begin{equation}\label{2.1.1}\lim_{n\rightarrow\infty}\sup_{y\in\mathbb{R}^N}\int_{B_r(y)}|u_n(x)|^2dx=0,\end{equation} for some $r>0$. Then $u_n\rightarrow0$ in
$L^q(\mathbb{R}^N)$ for any $q\in (2,2^*_\alpha)$.
\end{lemma}

For simplicity, in the sequel we call $\{u_n\}$ vanishing if $\{u_n\}$ satisfies (\ref{2.1.1}). We give some properties of $f$ in the following lemma, whose simple proof we omit.
\begin{lemma}\label{l2.1}

\noindent(i) If (f$_1$), (f$_3$) and (f$_4$) are satisfied, then $0<f(t)<l_0t$ for $t>0$.

\noindent(ii) If (f$_1$) and (f$_3$) are satisfied, then for $2<p< 2^*_\alpha$ and $\delta>0$, there exists $C_\delta=C_\delta(p)>0$ such that
\begin{equation}\label{2.2}|f(t)|\leq \delta|t|+C_\delta|t|^{p-1},\ |F(t)|\leq \delta|t|^2+C_\delta|t|^p,\ \forall t\in\mathbb{R}.\end{equation}

\noindent(iii) If (f$_1$) and (f$_4$) are satisfied, then for any $t>0$ \begin{equation}\label{2.3}tf'(t)>f(t),  \bar{F}(t):=\frac{1}{2}f(t)t-F(t)>0, \ \text{and} \ t\mapsto\bar{F}(t)\ \text{is  increasing}. \end{equation}

\noindent(iv) If (f$_3$)  is satisfied, then for any $u\in H^\alpha(\mathbb{R}^N)$\begin{equation}\label{2.4}
 \lim_{t\rightarrow\infty}\frac{1}{t^2}\int_{\mathbb{R}^N}F(tu)dx=\frac{l_0}{2}\int_{\mathbb{R}^N}u^2dx,\quad
 \lim_{t\rightarrow\infty}\frac{1}{t}\int_{\mathbb{R}^N}f(tu)udx=l_0\int_{\mathbb{R}^N}u^2dx.\end{equation} \end{lemma}

Now we describe the variational
framework.  The Nehari manifold $M_\epsilon$
of $I_\epsilon$ is
$$ M_\epsilon=\{u\in
H^{\alpha}(\mathbb{R}^N)\backslash\{0\}: \langle I'_\epsilon(u),u\rangle=0\}, $$
 and the least energy on $M_\epsilon$ is defined by $$c_\epsilon:=\inf_{M_\epsilon}I_\epsilon.$$
As mentioned in Section 1, we cannot use the method of Nehari manifold directly. We shall utilize a restricted set $\Theta_\epsilon$ given in (\ref{1.3}). Firstly, we  assert  $\Theta_\epsilon\neq\emptyset.$
Indeed, choose $u\in C^\infty_0(\mathbb{R}^N)$ such that $|u|_2=1$. Set $u_t(x):=t^{\frac{N}{2}}u(tx)$, $t>0$. Then $|u_t|_2=1$ and $[ u_t]_\alpha=t^\alpha[u]_\alpha\rightarrow0$ as $t\rightarrow0$. Hence, by (f$_3$) there exists small $t_0>0$ such that $$[ u_{t_0}]^2_\alpha+\int_{\mathbb{R}^N}V(\epsilon x)u^2_{t_0}(x)dx-l_0|u_{t_0}|^2_2\leq[ u_{t_0}]^2_\alpha+|V|_\infty-l_0<0.$$ Then $\Theta_\epsilon\neq\emptyset$.

\begin{lemma} \label{l2.2} Let (V$_1$) and (f$_1$)-(f$_4$) hold. For any $\epsilon>0$ the following results hold:

\noindent(i) $M_\epsilon\subset\Theta_\epsilon$ and  $M_\epsilon$ is bounded away from $0$;

\noindent(ii) for any $u\in \Theta_\epsilon$, set $\zeta_\epsilon(t):=I_\epsilon(tu)$. Then there exists a
unique $t_\epsilon:=t_\epsilon(u)>0$ such that $t_\epsilon u\in M_\epsilon$ and $I_\epsilon(t_\epsilon u)=\max_{t>0}I_\epsilon(t u)$.
Moreover,  \begin{equation}\label{2.5}c_\epsilon=\inf_{u\in
M_\epsilon}I_\epsilon(u)=\inf_{u\in \Theta_\epsilon}\max_{t>0}I_\epsilon(tu);\end{equation}

\noindent(iii) $c_\epsilon\geq r>0$.

\end{lemma}
{\bf Proof}: (i) Using Lemma \ref{l2.1} (i), we have $M_\epsilon\subset\Theta_\epsilon$. By (\ref{2.2}) we know for some $p\in(2,2^*_\alpha)$,
$$\|u\|^2_\epsilon\leq\delta\|u\|^2_\epsilon+C_\delta\|u\|^p_\epsilon,\ \ \forall u\in M_\epsilon.$$
Let $\delta=\frac12$, there holds $\|u\|_\epsilon\geq\varrho>0$.

(ii) From (\ref{2.2}) we have $\zeta_\epsilon(t)>0$ when $t$ is
sufficiently small.  In view of (\ref{2.4}) we obtain
 $$\zeta_\epsilon(t)=\frac{t^2}{2}\Bigl(\|u\|^2_\epsilon-l_0|u|^2_2+o_t(1)\Bigr),$$
as $t\rightarrow\infty$. Since $u\in \Theta_\epsilon$, we have $\zeta_\epsilon(t)<0$ for large $t$. Then $\zeta_\epsilon$ has maximum points in $(0,\infty)$ and (f$_4$) implies that the maximum point is unique. Therefore, (\ref{2.5}) follows.

(iii) Using (\ref{2.2}), there exist $\rho, r>0$ such that $\inf_{S_\rho}I_\epsilon\geq r$, where $S_\rho:=\{u\in H^\alpha(\mathbb{R}^N):\|u\|_\epsilon=\rho\}$.
For any $u\in M_\epsilon$, there is $t>0$ such that $tu\in S_\rho$. Using the conclusion (ii), we infer $I_\epsilon(u)\geq I_\epsilon(tu)$, then $\inf_{S_\rho}I_\epsilon\leq\inf_{M_\epsilon}I_\epsilon$. Hence $c_\epsilon\geq r>0$.
\ \ \ \ $\Box$
\begin{lemma} \label{l2.2.1} Under the assumptions of Lemma \ref{l2.2}, the following results hold:

\noindent(i) $M_\epsilon$ is a $C^1$ submanifold of $H^\alpha(\mathbb{R}^N)$;

\noindent(ii) if $u$ is a critical point of $I_\epsilon$ constrained on $M_\epsilon$, then $u$ is
a critical point of $I_\epsilon$.
\end{lemma}
{\bf Proof}: (i) For any $u\in H^\alpha(\mathbb{R}^N)$, denote \begin{equation}\label{2.7}J_\epsilon(u):=\langle I'_\epsilon(u),u\rangle.\end{equation} Observe that $$\langle J'_\epsilon(u),u\rangle=
2[u]^2_\alpha+2\int_{\mathbb{R}^N}V(\epsilon x)u^2dx-\int_{\mathbb{R}^N} [f(u)u+f'(u)u^2]dx.
$$
Then for each $u\in M_\epsilon
$, from (\ref{2.3}) we derive
\begin{equation}\label{2.6}\langle J'_\epsilon(u),u\rangle=\langle J'_\epsilon(u),u\rangle-2J_\epsilon(u)=\int_{\mathbb{R}^N}[f(u)u-f'(u)u^2]dx<0.\end{equation}
Thus, $ M_\epsilon$ is a $C^1$ submanifold of $H^\alpha(\mathbb{R}^N)$.

(ii)\ If $u$ is a critical point of $I_\epsilon$ constrained on $M_\epsilon
$, then there exists $\lambda\in \mathbb{R}$ such that $I'_\epsilon(u)=\lambda J'_\epsilon(u)$. Consequently
$0=\langle I'_\epsilon(u),u\rangle=\lambda\langle J'_\epsilon(u),u\rangle.$
By (\ref{2.6}) we get $\lambda=0$, and so $I'_\epsilon(u)=0$.\ \ \ \ $\Box$

\begin{lemma}\label{l2.3} Let (V$_1$) and (f$_1$)-(f$_5$) hold. Assume that $\{u_n\}\subset M_\epsilon
$ satisfies $I_\epsilon(u_n)\leq d$ with $d\in (0,\infty)$. Then $\{u_n\}$ is bounded in $H^\alpha(\mathbb{R}^N)$ and nonvanishing.
\end{lemma}
{\bf Proof}: Argue by contradiction we may assume  $\|u_n\|_\epsilon\rightarrow\infty$. Set $v_n:=\frac{u_n}{\|u_n\|_\epsilon}$. For any $t>0$, from (\ref{2.2}) and Lemma \ref{l2.2} (ii) we infer
\begin{equation}\label{2.7.0}d\geq I_\epsilon(u_n)\geq I_\epsilon(tv_n)\geq\frac{t^2}{2}-\int_{\mathbb{R}^N}F(tv_n)dx\geq\frac{t^2}{4}-
C t^p|v_n|^p_p,\end{equation}
where $p\in(2,2^*_\alpha)$.
If $\{v_n\}$ is vanishing, by Lemma \ref{l2.0.2} we get $|v_n|_p\rightarrow0$. There is a contradiction of (\ref{2.7.0}) if we choose $t>\sqrt{4d+1}$. If $\{v_n\}$ is nonvanishing, then there exist $y_n\in\mathbb{R}^N$ and $\delta_0>0$ such that $\int_{B_1(y_n)}v^2_ndx\geq\delta_0$. Set $\tilde{v}_n(\cdot)=\tilde{v}_n(\cdot+y_n)$. Then $\tilde{v}_n\rightharpoonup \tilde{v}\neq0$ and  we assume that $\tilde{v}(x)\neq0$, $x\in \Omega$, where $\Omega$ is a subset of $\mathbb{R}^N$ with positive measure. Denote $\tilde{u}_n=\|u_n\|_\epsilon \tilde{v}_n$. There holds $|\tilde{u}_n(x)|\rightarrow\infty$, $x\in \Omega$. According to Fatou lemma, (\ref{2.3}) and (f$_5$), we deduce
$$\aligned d&\geq I_\epsilon(u_n)-\frac12\langle I'_\epsilon(u_n),u_n\rangle=\int_{\mathbb{R}^N}\bigl[\frac12f(u_n)u_n-F(u_n)\bigl]dx\\&
=\int_{\mathbb{R}^N}\bigl[\frac12f(\tilde{u}_n)\tilde{u}_n-F(\tilde{u}_n)\bigl]dx=+\infty.\endaligned$$
This is a contradiction. Therefore, $\{u_n\}$ is bounded in $H^\alpha(\mathbb{R}^N)$. Next we show $\{u_n\}$ is nonvanishing.  If $\{u_n\}$ is vanishing, by Lemma \ref{l2.2} we conclude $u_n\rightarrow0$ in $L^p(\mathbb{R}^N)$ with $p\in(2,2^*_\alpha)$. Using (\ref{2.2}) and the boundedness of $\|u_n\|_\epsilon$  we have $\int_{\mathbb{R}^N}f(u_n)u_ndx\rightarrow0$. Then $\langle I'_\epsilon(u_n),u_n\rangle=0$ implies that $\|u_n\|_\epsilon\rightarrow0$, and so $I_\epsilon(u_n)\rightarrow0$. This is impossible since $I_\epsilon(u_n)\geq c_\epsilon\geq r$ using Lemma \ref{l2.2} (iii).  \ \ \ \
$\Box$

\begin{lemma}\label{l2.4} Let (V$_1$) and (f$_1$)-(f$_5$) hold. If $\{u_n\}$ is a PS sequence of $I_\epsilon$ constrained on $M_\epsilon
$, then $\{u_n\}$ is
a PS sequence of $I_\epsilon$.   \end{lemma}
{\bf Proof}: By Lemma \ref{l2.3}, $\{u_n\}$ is bounded in $H^{\alpha}(\mathbb{R}^N)$ and nonvanishing.
 Then there exists $y_n\in\mathbb{R}^N$ such that $\tilde{u}_n(\cdot)=u_n(\cdot+y_n)\rightharpoonup u\neq0$. So there exists $\Omega_0\subset\mathbb{R}^N$ with positive measure such that $|u(x)|>0$, $x\in \Omega_0$. Since $\{u_n\}$ is a PS sequence of $I_\epsilon$ constrained on $M_\epsilon
$, there exists $\sigma_n\in\mathbb{R}$ such that
$$o_n(1)=\nabla(I_\epsilon|_{M_\epsilon
})(u_n)=I'_\epsilon(u_n)-\sigma_nJ'_\epsilon(u_n),$$ where $J_\epsilon$ is given in (\ref{2.7}). Then
$$o_n(1)=\langle I'_\epsilon(u_n),u_n\rangle-\sigma_n\langle J'_\epsilon(u_n),u_n\rangle=-\sigma_n\langle J'_\epsilon(u_n),u_n\rangle.$$
From Fatou lemma, (\ref{2.3}) and (\ref{2.6}) it follows that
\begin{equation}\label{3.1}\aligned|\langle J'_\epsilon(u_n),u_n\rangle|
&=\Bigl|\int_{\mathbb{R}^N}[f(u_n)u_n-f'(u_n)u^2_n]dx\Bigr|
=\Bigl|\int_{\mathbb{R}^N}[f(\tilde{u}_n)\tilde{u}_n-f'(\tilde{u}_n)\tilde{u}^2_n]dx\Bigr|\\
 &\geq\int_{\Omega_0}[f'(u)u^2-f(u)u]dx>0.\endaligned\end{equation}
Then $\sigma_n\rightarrow0$. Moreover, in view of the boundedness of $\|u_n\|_\epsilon$, (f$_2$) and (\ref{2.2}) we have
$$\Bigl|2\langle u_n,\phi\rangle_\epsilon-\int_{\mathbb{R}^N}\bigl[f'(u_n)u_n+f(u_n)\bigr]\phi dx\Bigr|\leq C(\|u_n\|_\epsilon+\|u_n\|^{p-1}_\epsilon+\|u_n\|^{q-1}_\epsilon)\|\phi\|\leq C\|\phi\|.$$
 Then $\{J'_\epsilon(u_n)\}$ is bounded. Hence $I'_\epsilon(u_n)\rightarrow0$.\ \ \ \ \ $\Box$
\section{The autonomous problem}

In this section we are concerned with the autonomous equation. Precisely, for
any constant $a$ with $0<a<l_0$, consider
\begin{equation}\label{3.1}
(-\Delta)^\alpha u+au=f(u),\ \ \ x\in \mathbb{R}^{N}.
\end{equation}
The functional of (\ref{3.1}) is denoted by
$$I_{a}(u)=\frac{1}{2}([u]^{2}_\alpha+a|u|^{2}_2)-\int_{\mathbb{R}^{N}}F(u)dx.$$
The Nehari manifold corresponding to (\ref{3.1}) is defined by
$$ M_{a}=\{u\in
H^{\alpha}(\mathbb{R}^N)\backslash\{0\}:\langle I'_{a}(u),u\rangle=0\}, $$
 and the least energy on $M_{a}$ is defined by $c_{a}:=\inf_{M_{a}}I_{a}.$ To give a new characterization of $c_a$, set
$$\Theta_{a}:=\{u\in H^\alpha(\mathbb{R}^N): [u]^2_\alpha+a|u|^2_2-l_0|u|^2_2<0\}.$$
Replacing $I_\epsilon$, $M_\epsilon
$, $c_\epsilon$ and $\Theta_\epsilon$ with $I_a$, $M_a$, $c_a$ and $\Theta_a$ respectively, it is easy to see that Lemmas \ref{l2.2}-\ref{l2.4} still hold and below we state some of them.
\begin{lemma} \label{l3.1} For any constant $a\in(0,l_0)$, the following results hold:

 \noindent(i) $M_a\subset \Theta_a$ and $M_a$ is bounded away from $0$;

 \noindent(ii)  for all $u\in \Theta_a$, there exists a
unique $t_{u}>0$ such that $t_{u}u\in M_a$ and
$I_a(t_{u}u)=\max_{t>0}I_a(tu)$. Moreover,
\begin{equation}\label{3.2}c_{a}=\inf_{u\in M_{a}}I_{a}(u)=\inf_{w\in \Theta_a}\max_{t>0}I_{a}(tw)>0;\end{equation}

 \noindent(iii) If $\{u_n\}$ is a PS sequence of $I_a$ constrained on $M_a
$, then $\{u_n\}$ is bounded in $H^\alpha(\mathbb{R}^N)$ and nonvanishing. Moreover, $\{u_n\}$ is
a PS sequence of $I_a$;

\noindent(iv) if $u$ is a critical point of $I_a$ constrained on $M_a$, then $u$ is
a critical point of $I_a$.
\end{lemma}

By \cite{CW}, we know the problem (\ref{3.1}) has a positive ground state $u_{a}$. We would like to point out that, the result in  \cite{CW} was obtained under the assumption that $f$ is odd in $\mathbb{R}$. Since we just consider positive solutions, extending the nonlinearity $f$ to an odd
function on $\mathbb{R}$, we can apply \cite{CW} to our setting. In this paper, in order to investigate the behavior of the least energy for different parameter $a$, we give another proof using the method of Nehari manifold and show that $c_a$ is attained.

\begin{lemma}\label{l3.2} For any $0<a<l_0$, problem (\ref{3.1}) has a positive ground state $u$ satisfying $I_a(u)=c_a$.
\end{lemma}
{\bf Proof}: Assume that $u_n\in M_a$ satisfies that $I_a(u_n)\rightarrow
c_a$. By the Ekeland variational principle, we may
suppose that $\nabla(I_a|_{M_a})(u_n)\rightarrow0$. From
Lemma \ref{l3.1} (iii) it follows that  $\{u_n\}$ is bounded in $H^\alpha(\mathbb{R}^N)$ and  $I'_a(u_n)\rightarrow 0$.  In addition,
$\{u_n\}$ is nonvanishing.
Then there exist $x_n\in
\mathbb{R}^N$ and $\delta_0>0$ such that
\begin{equation}\label{3.4} \int_{B_1(x_n)}u^2_n(x)dx>\delta_0.\end{equation}
Since $I_a$ and $M_a$ are invariant by the translation of the form $v\mapsto v(\cdot-z)$, $z\in\mathbb{R}^N$.
Without loss of generality we suppose $\{x_n\}$ is bounded. By (\ref{3.4}) we may suppose $u_n\rightharpoonup {u_0}\neq0$ in $H^{\alpha}(\mathbb{R}^N)$, $u_n\rightarrow {u_0}$ in $L^{2}_{loc}(\mathbb{R}^N)$ and $u_n(x)\rightarrow{u_0}(x)$ a.e. in $\mathbb{R}^N$. Then
\begin{equation}\label{3.5}\aligned c_a\leq& I_a(u_0)-\frac12\langle I'_a(u_0),u_0\rangle=\int_{\mathbb{R}^N}[\frac12 f(u_0)u_0-F(u_0)]dx\\ \leq &\liminf_{n\rightarrow\infty}\Bigl\{\int_{\mathbb{R}^N}[\frac12 f(u_n)u_n-F(u_n)]dx\Bigr\}=\liminf_{n\rightarrow\infty}\Bigl\{I_a(u_n)-\frac12\langle I'_a(u_n),u_n\rangle\Bigr\}=c_a.\endaligned\end{equation}
So $I_a(u_0)=c_a$ and $u_0$ is a ground state of (\ref{3.1}).
Using a  standard argument, we may assume that $u_0>0$. This ends the proof. \ \ \ \ $\Box$

By Lemma \ref{l3.2} and the characterization (\ref{3.2}), it is easy to have the following result whose proof we omit.
\begin{lemma}\label{l3.3}Assume $0<a_1<a_2<l_0$. Then
$c_{a_1}<c_{a_2}$.
\end{lemma}

\section{Proof of Theorem 1.1}

In this section, we are devoted to proving Theorem 1.1 and we always assume that (V$_1$) and (f$_1$)-(f$_5$) are satisfied.

 Let $w$ be a positive ground state of problem (\ref{3.1}) with $a=V_0$ and $\eta$ be a smooth non-increasing
function defined in $[0,+\infty)$ such that $\eta(s)=1$ if $0\leq s\leq\frac12$ and $\eta(s)=0$ if $s\geq1$. For
any $y\in \mathcal{V}$, we define
$$\psi_{\epsilon,y}(x)=\eta(|\epsilon x-y|)w\bigl(\frac{\epsilon x-y}{\epsilon}\bigr).$$
We claim that $\psi_{\epsilon,y}\in \Theta_\epsilon$ for small enough $\epsilon$. Indeed, the Lebesgue's convergence theorem implies
\begin{equation}\label{4.2.1}\aligned&\lim_{\epsilon\rightarrow0}[ \psi_{\epsilon,y}]^2_\alpha= [w]^2_\alpha, \ \lim_{\epsilon\rightarrow0}\int_{\mathbb{R}^N}V(\epsilon x)|\psi_{\epsilon,y}|^2dx= V_0|w|^2_2,\\ &\lim_{\epsilon\rightarrow0}\int_{\mathbb{R}^N}f(\psi_{\epsilon,y})\psi_{\epsilon,y}dx=\int_{\mathbb{R}^N}f(w)wdx.\endaligned\end{equation}
Then
\begin{equation}\label{4.7.1}\aligned&[ \psi_{\epsilon,y}]^2_\alpha+\int_{\mathbb{R}^N}V(\epsilon x)|\psi_{\epsilon,y}|^2dx-l_0|\psi_{\epsilon,y}|^2_2\\=&[ w]^2_\alpha+V_0|w|^2_2-l_0|w|^2_2+o_\epsilon(1)=\int_{\mathbb{R}^N}f(w)wdx-l_0|w|^2_2+o_\epsilon(1).\endaligned\end{equation}
By Lemma \ref{l2.1} (i) there exists $\epsilon(y)>0$ such that for any $\epsilon<\epsilon(y)$,
$\psi_{\epsilon,y}\in \Theta_\epsilon$. Since $\mathcal{V}$ is compact, there exists $\epsilon_1>0$ independent on $y$ such that for each $\epsilon<\epsilon_1$, $\psi_{\epsilon,y}\in \Theta_\epsilon$. In this section, we always assume $\epsilon<\epsilon_1$.
From Lemma \ref{l2.2} (ii) we know there exists $t_\epsilon>0$ such that $\max_{t\geq0}{I}_\epsilon(t\psi_{\epsilon,y})={I}_\epsilon(t_\epsilon\psi_{\epsilon,y})$.
Define $\gamma_{\epsilon}: \mathcal{V}\rightarrow M_\epsilon$ by $\gamma_{\epsilon}(y)=t_\epsilon\psi_{\epsilon,y}$.
By the construction, $\gamma_{\epsilon}(y)$ has a compact support for any $y\in \mathcal{V}$.

\begin{lemma}\label{l4.1}
$\lim_{\epsilon\rightarrow0}{I}_\epsilon(\gamma_{\epsilon}(y))={c}_{V_0}$ uniformly in $y\in \mathcal{V}$.
\end{lemma}
{\bf Proof}: Suppose by contradiction that there exist $\delta_0>0$, $y_n\in \mathcal{V}$ and $\epsilon_n\rightarrow0$ such that
\begin{equation}\label{4.2.3}|{I}_{\epsilon_n}(\gamma_{\epsilon_n}(y_n))-{c}_{V_{0}}|\geq \delta_0.\end{equation}
Now we claim that \begin{equation}\label{4.9.0}\lim_{n\rightarrow\infty}t_{\epsilon_n}=1.\end{equation} In fact, since $t_{\epsilon_n}\psi_{\epsilon_n,y_n}\in M_{\epsilon_n}$ and $M_{\epsilon_n}$ is bounded away from 0, $t_{\epsilon_n}$ cannot go zero.
Assume that $t_{\epsilon_n}\rightarrow A>0$. If $A\geq1+2\delta>1$ may be $+\infty$, then
from (f$_4$) and Lebesgue's convergence theorem it follows that \begin{equation}\label{4.2.2}\aligned\int_{\mathbb{R}^N}\frac{f(t_{\epsilon_n} \psi_{\epsilon_n,y_n})\psi_{\epsilon_n,y_n}}{t_{\epsilon_n}}dx&\geq \int_{\mathbb{R}^N}\frac{f((1+\delta) \psi_{\epsilon_n,y_n})\psi_{\epsilon_n,y_n}}{1+\delta}dx
\\&=\int_{\mathbb{R}^N}\frac{f((1+\delta) w)w}{1+\delta}dx+o_n(1)>\int_{\mathbb{R}^N} f(w)wdx+o_n(1),\endaligned\end{equation}
for large $n$. However, since $t_{\epsilon_n}\psi_{\epsilon_n,y_n}\in M_{\epsilon_n}$ and (\ref{4.2.1}) we have
$${t^{-1}_{\epsilon_n}}\int_{\mathbb{R}^N}f(t_{\epsilon_n} \psi_{\epsilon_n,y_n})\psi_{\epsilon_n,y_n}dx=[w]^2_\alpha+V_0|w|^2_2+o_n(1),$$
which contradicts with (\ref{4.2.2}). Similarly, we can exclude the case that $0<A<1$. So $A=1$.
Then by (\ref{4.2.1}) we infer ${I}_{\epsilon_n}(\gamma_{\epsilon_n}(y_n))\rightarrow {I}_{V_0}(w)={c}_{V_0}$, contradicting with (\ref{4.2.3}).\ \ \ \ $\Box$

\begin{lemma}\label{l4.2}$\lim_{\epsilon\rightarrow0}c_\epsilon=c_{V_0}$.\end{lemma}
{\bf Proof}: By Lemma \ref{l4.1} we have that $\limsup_{\epsilon\rightarrow0}c_\epsilon\leq c_{V_0}$. Combining (\ref{2.5}) and (\ref{3.2}), we deduce $\liminf_{\epsilon\rightarrow0}c_\epsilon\geq c_{V_0}$. Then $\lim_{\epsilon\rightarrow0}c_\epsilon=c_{V_0}$.\ \ \ \ \ $\Box$

From (V$_1$), we know $V_0<V_\infty$. Choose $\xi\in (V_0,V_\infty)$. Lemma \ref{l3.3} implies that $c_{V_0}<c_{\xi}<c_{V_\infty}$. By Lemma \ref{l4.2}, we may suppose
\begin{equation}\label{4.1}c_{\epsilon}\leq c_{\xi}<c_{V_\infty},\end{equation}
for small $\epsilon$. Next we give a compactness lemma.
\begin{lemma}\label{l4.3} For any bounded (PS)$_c$ sequence $\{u_n\}$ for $I_\epsilon$ with $0<c\leq c_\xi<c_{V_\infty}$,  $\{u_n\}$ converges strongly in $H^\alpha(\mathbb{R}^N)$ after passing to a subsequence.
\end{lemma}
{\bf Proof}: Assume $u_n\rightharpoonup u$ in $H^\alpha(\mathbb{R}^N)$. Define $v_n=u_n-u$ and suppose that $v_n\not\rightarrow0$ in $H^\alpha(\mathbb{R}^N)$. As \cite[Lemma 2.14]{HZ} we have
\begin{equation}\label{4.2}\langle I'_\epsilon(v_n),v_n\rangle=\langle I'_\epsilon(u_n),u_n\rangle-\langle I'_\epsilon(u),u\rangle+o_n(1)=o_n(1),\end{equation}
\begin{equation}\label{4.3} I_\epsilon(v_n)=I_\epsilon(u_n)-I_\epsilon(u)+o_n(1).\end{equation}
If $\{v_n\}$ is vanishing, then by Lemma \ref{l2.0.2} we know $v_n\rightarrow0$ in $L^p(\mathbb{R}^N)$ with $p\in(2,2^*_\alpha)$. From (\ref{4.2}) we get $\|v_n\|_\epsilon\rightarrow0$, which contradicts with the hypothesis
$v_n\not\rightarrow0$ in $H^\alpha(\mathbb{R}^N)$. So $\{v_n\}$ is nonvanishing. Then there exists $y_n\in\mathbb{R}^N$ such that
 $\tilde{v}_n(\cdot)={v}_n(\cdot+y_n)\rightharpoonup \tilde{v}\neq0$ in $H^\alpha(\mathbb{R}^N)$.
Given $\varsigma>0$, from assumption (V$_1$), there exists $R=R(\varsigma)>0$ such that $V(\epsilon x)\geq V_\infty-\varsigma$, {for any} $|x|\geq R.$
Using the fact that $v_n\rightarrow0$ in $L^2(B_R(0))$,  we infer
\begin{equation}\label{4.5}\int_{\mathbb{R}^N}V(\epsilon x)v^2_ndx\geq V_\infty|v_n|^2_2-C\varsigma.\end{equation}
According to Lemma \ref{l2.1} (i), (\ref{4.2}) and Fatou Lemma we obtain
\begin{equation}\label{4.5.0}\aligned &l_0|v_n|^2_2-[ v_n]^2_\alpha-V_\infty|v_n|^2_2\geq l_0|v_n|^2_2-[ v_n]^2_\alpha-\int_{\mathbb{R}^N}V(\epsilon x)v^2_ndx-C\varsigma\\
=&l_0|v_n|^2_2-\int_{\mathbb{R}^N}f(v_n)v_ndx+o_n(1)-C\varsigma
=l_0|\tilde{v}_n|^2_2-\int_{\mathbb{R}^N}f(\tilde{v}_n)\tilde{v}_ndx+o_n(1)-C\varsigma\\
\geq& l_0|\tilde{v}|^2_2-\int_{\mathbb{R}^N}f(\tilde{v})\tilde{v}dx+o_n(1)-C\varsigma.\endaligned\end{equation}
In view of the arbitrariness of $\varsigma>0$ and $\tilde{v}\neq0$ we know $v_n\in \Theta_{V_\infty}$.
Let $t_n>0$ be such that $t_nv_n\in M_{V_\infty}$.

Case 1:
$\limsup_{n\rightarrow\infty}t_n>1.$ Up to a subsequence,  we assume that
$t_n\geq 1+\delta>1$ for all $n\in\mathbb{N}.$
By (\ref{4.2}), (\ref{4.5}) and $t_nv_n\in M_{V_\infty}$ we deduce
\begin{equation}\label{4.6}\int_{\mathbb{R}^N}
\bigl[\frac{f(t_nv_n)v_n}{t_n}-f(v_n)v_n\bigr]dx+o_n(1)=\int_{\mathbb{R}^N}(V_\infty-V(\epsilon x))|v_n|^2dx\leq C\varsigma+o_n(1).\end{equation}
 Since $\tilde{v}\neq0$, there exists a subset $\Omega\subset\mathbb{R}^N$ with positive measure such that $|\tilde{v}(x)|>0$ for all $x\in\Omega$.
Letting $n\rightarrow\infty$ in (\ref{4.6}), by (f$_4$) there holds
$$0<\int_{\Omega}
\bigl[\frac{f((1+\delta)\tilde{v}(x))\tilde{v}(x)}{1+\delta}
-f(\tilde{v}(x))\tilde{v}(x)\bigr]dx\leq C\varsigma,$$
which is absurd using the arbitrariness of $\varsigma$.

Case 2: $\limsup_{n\rightarrow\infty}t_n=1$.
After passing to a subsequence, we assume that $t_n\rightarrow1$.
By (\ref{4.5})  we infer
$$I_\epsilon(v_n)-I_{V_\infty}(t_nv_n)=I_\epsilon(v_n)-I_{V_\infty}(v_n)+o_n(1)\geq o_n(1)-C\varsigma. $$
From (\ref{4.3})
we have
$I_\epsilon(u_n)\geq c_{V_\infty}+o_n(1)-C\varsigma+I_\epsilon(u).$
Using (\ref{2.3}) we get
\begin{equation}\label{4.0.0}I_\epsilon(u)=I_\epsilon(u)-\frac{1}{2}\langle I'_\epsilon(u),u\rangle=\int_{\mathbb{R}^N}\bigl[\frac12f(u)u-F(u)\bigr]dx\geq0.\end{equation}
Then $c\geq c_{V_\infty}$, which contradicts with $c<c_{V_\infty}$.

Case 3: $\limsup_{n\rightarrow\infty}t_n=t_0<1$ and we may suppose  $t_n\rightarrow t_0<1$.
From (\ref{2.3}) we derive
$$\aligned c_{V_\infty}
&\leq I_{V_\infty}(t_nv_n)-\frac12\langle I'_{V_\infty}(t_nv_n), t_nv_n\rangle=\int_{\mathbb{R}^N}\bigl[\frac12 f(t_nv_n(x))t_nv_n(x)-F(t_nv_n(x))\bigr]dx\\&
\leq\int_{\mathbb{R}^N}\bigl[\frac12 f(v_n(x))v_n(x)-F(v_n(x))\bigr]dx
=I_\epsilon(v_n)-\frac1{2}\langle I'_\epsilon(v_n),v_n\rangle=c-I_\epsilon(u)+o_n(1).\endaligned$$
Hence,  $c-I_\epsilon(u)\geq c_{V_\infty},$ contradicting with (\ref{4.0.0}) and $c<c_{V_\infty}$. \ \ \ $\Box$

Since $f$ is asymptotically linear at infinity, it is not easy to show PS condition holds for $I_{V_0}$.  We need to investigate PS sequence carefully.
\begin{lemma}\label{l4.4} Assume that $v_n\in M_{V_0}$ satisfies $I_{V_0}(v_n)\rightarrow{c_{V_0}}$ and $v_n\rightharpoonup v\neq0$ in $H^\alpha(\mathbb{R}^N)$. Then $v_n\rightarrow v$ in $H^\alpha(\mathbb{R}^N)$.
\end{lemma}
{\bf Proof}:  Since $v_n$ is a minimizing sequence of $I_{V_0}|_{M_{V_0}}$, Ekeland variational principle implies that
$\nabla(I_{V_0}|_{M_{V_0}})(v_n)\rightarrow0$. From Lemma \ref{l3.1} (iii) we know $I'_{V_0}(v_n)\rightarrow0$. Set $w_n:=v_n-v$. Similar to (\ref{4.2}) and (\ref{4.3}) we have
\begin{equation}\label{6.3.1}\langle I'_{V_0}(w_n),w_n\rangle=o_n(1),\ \ I_{V_0}(w_n)=I_{V_0}(v_n)-I_{V_0}(v)+o_n(1).\end{equation}
If $\|w_n\|\geq\rho>0$, by (\ref{6.3.1}) we conclude that $w_n$ is nonvanishing, and then we claim $w_n\in \Theta_{V_0}$.
Indeed, since $w_n$ is nonvanishing, there exists $y_n\in\mathbb{R}^N$ such that $w_n(\cdot+y_n)\rightharpoonup w\neq0$ in $H^\alpha(\mathbb{R}^N)$. As (\ref{4.5.0}) we get$$\aligned l_0|w_n|^2_2-[w_n]^2_\alpha-V_0|w_n|^2_2\geq l_0|w|^2_2-\int_{\mathbb{R}^N}f(w)wdx+o_n(1).\endaligned$$
Then $w_n\in \Theta_{V_0}$.
By Lemma \ref{l3.1} (ii), there exists $t_n>0$ such that $t_n w_n\in M_{V_0}$. Due to Lemma \ref{l3.1} (i), $t_n$ cannot go zero. We next claim  $t_n\rightarrow1$. In fact,
by $t_n w_n\in M_{V_0}$ and (\ref{6.3.1})
 we obtain
$$o_n(1)=\int_{\mathbb{R}^N}\Bigl[\frac{f(t_nw_n)w_n}{t_n}-{f(w_n)w_n}\Bigr]dx.$$
If $t_n \rightarrow A>1$, then
as Case 1 in the proof of Lemma \ref{l4.3}, there will be a contradiction. Similarly, we can exclude the case $0<A<1$. Then $A=1$.
Moreover, by the assumption $v\neq0$, we have $v\in M_{V_0}$. From (\ref{6.3.1}) and $t_n\rightarrow1$  it follows that
$ I_{V_0}(t_n w_n)\leq o_n(1),$
which contradicts with
$I_{V_0}(t_n w_n)\geq c_{V_0}$.
Thus $v_n\rightarrow v$ in $H^\alpha(\mathbb{R}^N)$.\ \ \ \ $\Box$

For any $\delta>0$, let $\rho=\rho(\delta)>0$ be such that $\mathcal{V}_\delta\subset B_\rho(0)$. Define
$\chi:=\mathbb{R}^N\rightarrow\mathbb{R}^N$ as $\chi(x)=x$ for $|x|\leq\rho$ and $\chi(x)=\frac{\rho x}{|x|}$
for $|x|\geq\rho$. Moreover, we define $\beta_\epsilon: {M}_\epsilon\rightarrow \mathbb{R}^N$ by
$$\beta_\epsilon(u):=\frac{\int_{\mathbb{R}^N}\chi(\epsilon x)u^2dx}{\int_{\mathbb{R}^N}u^2dx}.$$
Since $\mathcal{V}\subset B_\rho(0)$, by the definition of $\chi$ and Lebesgue's convergence theorem, we deduce
\begin{equation}\label{4.16.0}\aligned\beta_\epsilon(\gamma_\epsilon(y))&=y+\frac{\int_{\mathbb{R}^N}\bigl(\chi(\epsilon x+y)-y\bigr)|w(x)\eta(|\epsilon x|)|^2dx}{\int_{\mathbb{R}^N}|w(x)\eta(|\epsilon x|)|^2dx}=y+o_\epsilon(1),\endaligned\end{equation}
as $\epsilon\rightarrow0$, uniformly for $y\in \mathcal{V}$. So $\lim_{\epsilon\rightarrow0}\beta_\epsilon(\gamma_\epsilon(y))=y$ uniformly for $y\in \mathcal{V}$.
To go on, we need the following compact lemma.
\begin{lemma}\label{l4.5}Let $\epsilon_n\rightarrow0$ and $\{u_n\}\subset M_{\epsilon_n}$ be such that $I_{\epsilon_n}(u_n)\rightarrow c_{V_{0}}$. Then there exists $\{y_n\}\subset\mathbb{R}^N$ such that $v_n(x)=u_n(x+y_n)$ has a convergent subsequence in $H^\alpha(\mathbb{R}^N)$. Moreover, up to a subsequence, $\tilde{y}_n:=\epsilon_ny_n\rightarrow y_0\in\mathcal{V}$.
\end{lemma}
{\bf Proof}: As in the proof of Lemma \ref{l2.3}, $\{u_n\}$ is bounded in $H^\alpha(\mathbb{R}^N)$ and non-vanishing. Then
there exists $y_n\in\mathbb{R}^N$ such that
\begin{equation}\label{4.5.1}\int_{B_1(y_n)}u^2_n(x)dx>\delta.\end{equation}
We may assume that $v_n(\cdot)=u_n(\cdot+y_n)\rightharpoonup v\neq0$ in
$H^\alpha(\mathbb{R}^N)$. Since $u_n\in M_{\epsilon_n}\subset \Theta_{\epsilon_n}\subset\Theta_{V_{0}}$, there exists $t_n>0$ such that $\tilde{v}_n=t_nv_n\in M_{V_{0}}$.
Then
$ c_{V_0}\leq I_{V_{0}}(\tilde{v}_n)
\leq I_{\epsilon_n}(t_nu_n)\leq I_{\epsilon_n}(u_n).$
Hence, $\lim_{n\rightarrow\infty} I_{V_0}(\tilde{v}_n)=c_{V_0}$. Since $v_n\rightharpoonup v\neq0$ in
$H^\alpha(\mathbb{R}^N)$, there exists $C>0$ such that $\|v_n\|_{\epsilon_n}>C$. Therefore, $\|\tilde{v}_n\|_{\epsilon_n}=t_n\|v_n\|_{\epsilon_n}>t_nC$. By Lemma \ref{l3.1} (iii) we infer that $\{\tilde{v}_n\}$ is bounded in $H^\alpha(\mathbb{R}^N)$. So $\{t_n\}$ is bounded. Suppose  $t_n\rightarrow t\geq0$. Since $M_{V_0}$ is bounded away from $0$, from $\tilde{v}_n\in M_{V_0}$ we deduce $t>0$.  Let $\tilde{v}_n\rightharpoonup\tilde{v}$ in $H^\alpha(\mathbb{R}^N)$. Then $\tilde{v}=tv\neq0$.
Using Lemma \ref{l4.4} we get $\tilde{v}_n\rightarrow \tilde{v}$ in $H^\alpha(\mathbb{R}^N)$.

Set $\tilde{y}_n:=\epsilon_ny_n$. Next we claim that $\{\tilde{y}_n\}$ is bounded. Indeed, suppose by contradiction that $|\tilde{y}_n|\rightarrow\infty$. From $\tilde{v}_n\rightarrow\tilde{v}$ in $H^\alpha(\mathbb{R}^N)$ and ${V_0}<V_{\infty}$ we conclude
$$\aligned
I_{V_{\infty}}(\tilde{v})&=\lim_{n\rightarrow\infty}I_{V_{\infty}}(\tilde{v}_n)
\leq\liminf_{n\rightarrow\infty}\Bigl\{\frac12[\tilde{v}_n]^2_\alpha+\frac12\int_{\mathbb{R}^{N}}V(\epsilon_n x+\tilde{y}_n){\tilde{v}}^2_ndx-\int_{\mathbb{R}^{N}}F(\tilde{v}_n(x))dx\Bigr\}\\&
=\liminf_{n\rightarrow\infty}I_{\epsilon_n}(t_nu_n)
\leq\liminf_{n\rightarrow\infty}I_{\epsilon_n}(u_n)=c_{V_0}.\endaligned$$
This is impossible since $c_{V_0}=I_{V_0}(\tilde{v})<I_{V_{\infty}}(\tilde{v})$. So $\{\tilde{y}_n\}$ is bounded, and up to a subsequence, suppose that $\tilde{y}_n\rightarrow y_0$ in $\mathbb{R}^N$. If $y_0\not\in \mathcal{V}$, then $V(y)>V_0$ and we obtain a contradiction by the same arguments made above. So $y_0\in \mathcal{V}$.\ \ \ \ $\Box$

Let $\omega(\epsilon)$ be any positive function satisfying $\omega(\epsilon)\rightarrow0$ as $\epsilon\rightarrow0$. Define the set $$\Sigma_\epsilon:=\{u\in M_\epsilon: I_\epsilon(u)\leq c_{V_0}+\omega(\epsilon)\}.$$
For any $y\in\mathcal{V}$, we deduce from Lemma \ref{l4.1} that $\omega(\epsilon)=|I_\epsilon(\gamma_\epsilon(y))-c_{V_0}|\rightarrow0$ as $\epsilon\rightarrow0$. Thus $\gamma_\epsilon(y)\in \Sigma_\epsilon$ and $\Sigma_\epsilon\neq\emptyset$ for small $\epsilon>0$.

\begin{lemma}\label{l4.6} For any $\delta>0$,  $\lim_{\epsilon\rightarrow0}\sup_{u\in\sum_\epsilon}dist(\beta_\epsilon(u),\mathcal{V}_\delta)=0$.
\end{lemma}
{\bf Proof}: Set $\epsilon_n\rightarrow0$. By definition, there
exists $\{u_n\}\subset\Sigma_{\epsilon_n}$ such that $$dist(\beta_{\epsilon_n}(u_n),\mathcal{V}_\delta)
=\sup_{u\in\Sigma_{\epsilon_n}}dist(\beta_{\epsilon_n}(u),\mathcal{V}_\delta)+o_n(1).$$ It suffices to find a sequence $\{\tilde{y}_n\}\subset\mathcal{V}_\delta$ satisfying $|\beta_{\epsilon_n}(u_n)-\tilde{y}_n|=o_n(1)$. Since $\{u_n\}\subset\Sigma_{\epsilon_n}$, it is easy to see that
$I_{\epsilon_n}(u_n)\rightarrow c_{V_0}$. By Lemma \ref{l4.5}, there exists $\{y_n\}\subset\mathbb{R}^N$ such that $\tilde{y}_n=\epsilon_n y_n\in \mathcal{V}_\delta$ for $n$ sufficiently large.
Hence
$$\beta_{\epsilon_n}(u_n)
=\tilde{y}_n+\frac{\int_{\mathbb{R}^N}(\chi(\epsilon_nz+\tilde{y}_n)-\tilde{y}_n)u^2_n(z+{y}_n)dz}
{\int_{\mathbb{R}^N}u^2_n(z+{y}_n)dz}.$$
Since $\epsilon_nz+\tilde{y}_n\rightarrow y_0\in\mathcal{V}$, we have $\beta_{\epsilon_n}(u_n)
=\tilde{y}_n+o_n(1)$ and $\{\tilde{y}_n\}$ is what we need. \ \ \ \ $\Box$

Now we are ready to show the multiplicity of positive solutions.
\begin{lemma}\label{l4.7}The equation (\ref{2.1}) has at
least $cat_{\mathcal{V}_\delta}(\mathcal{V})$ positive solutions for small $\epsilon>0$.\end{lemma}
{\bf Proof}:
For fixed $\delta>0$,
from Lemmas \ref{l4.1} and \ref{l4.6}, we know that there exists $\epsilon_\delta>0$ such that for any $\epsilon\in (0,\epsilon_\delta)$, the diagram
$\mathcal{V}\xrightarrow{\gamma_\epsilon}\Sigma_{\epsilon}\xrightarrow{\beta_\epsilon}\mathcal{V}_\delta,$
is well defined. By (\ref{4.16.0}), we see that
$\lim_{\epsilon\rightarrow0}\beta_\epsilon(\gamma_\epsilon(y))=y$ {uniformly in}\ $y\in\mathcal{V}.$
For $\epsilon>0$ small enough, we denote $\beta_\epsilon(\gamma_\epsilon(y))=y+\nu(y)$ for $y\in \mathcal{V}$,
where $|\nu(y)|<\frac{\delta}{2}$ uniformly for $y\in \mathcal{V}$. Define $H(t,y)=y+(1-t)\nu(y)$. Then
$H:[0,1]\times\mathcal{V}\rightarrow\mathcal{V}_\delta$ is continuous. Obviously, $H(0,y)=\beta_\epsilon(\gamma_\epsilon(y))$, $H(1,y)=y$ for all $y\in \mathcal{V}$.
Thus we obtain that the composite mapping $\beta_{\epsilon}\circ \gamma_{\epsilon}$ is homotopic to the inclusion map $id: \mathcal{V}\rightarrow\mathcal{V}_\delta$. Using  \cite[Lemma 2.2]{CL} there holds
\begin{equation}\label{5.24}cat_{\Sigma_{\epsilon}}(\Sigma_{\epsilon})\geq cat_{\mathcal{V}_\delta}(\mathcal{V}).\end{equation}

Choose a function $\omega(\epsilon)>0$ such that $\omega(\epsilon)\rightarrow0$
as $\epsilon\rightarrow0$ and $c_{V_0}+\omega(\epsilon)$ is not a critical level for $I_\epsilon$. From the definition of $\Sigma_{\epsilon}$ and $c_{V_0}<c_{V_\infty}$, Lemmas \ref{l4.3} and \ref{l2.4} imply that  $I_\epsilon$ satisfies the Palais-Smale condition on
$\Sigma_{\epsilon}$. It follows from standard Ljusternik-Schnirelmann category theory that $I_\epsilon$ has at least $cat_{\Sigma_{\epsilon}}(\Sigma_{\epsilon})$ critical points on $\Sigma_{\epsilon}$. In view of  Lemma \ref{l2.2.1} (ii), (\ref{5.24}) and (f$_1$), $I_\epsilon$ has at least $cat_{\mathcal{V}_\delta}(\mathcal{V})$ positive critical points.
\ \ \ \ $\Box$

Next we study the concentration of solutions obtained in Lemma \ref{l4.7}.
Set $\epsilon_n\rightarrow0$. Let $u_n$ be the positive solutions of (\ref{2.1}) obtained in Lemma \ref{l4.7}. Then  $v_n(x):=u_n(x+y_n)$, where $y_n$ is given in Lemma \ref{l4.5}, are the associated solutions of the  equation
\begin{equation}\label{4.7}(-\Delta)^\alpha u+V_n(x)u=f(u), \ \ x\in\mathbb{R}^N,\end{equation}
where $V_n(x)=V(\epsilon_n x+\epsilon_n y_n)$.

\begin{lemma}\label{l5.9}
There exist $C',\delta'>0$ such that $C'\geq|{v}_n|_\infty\geq\delta'$ for all $n\in\mathbb{N}$, and \begin{equation}\label{4.18.0}\lim_{|x|\rightarrow\infty}{v}_n(x)=0\ \ \text{uniformly in }n\in\mathbb{N}.\end{equation}\end{lemma}
{\bf Proof}: As \cite[Lemma 5.1]{Ambro} or \cite[Lemma 4.1]{HZ1}, we have
\begin{equation}\label{5.9.1}|v_n|_\infty\leq C'.\end{equation}
In view of Lemma \ref{l4.5}, we know $v_n\rightarrow v\neq0$ in $H^\alpha(\mathbb{R}^N)$ and $\epsilon_ny_n\rightarrow y_0\in \mathcal{V}$.
By (\ref{4.5.1}), we deduce
$|v_n|_\infty\geq \delta'.$
Below we show (\ref{4.18.0}).
Since $v_n\rightarrow v$ in $H^\alpha(\mathbb{R}^N)$, from (\ref{5.9.1}) and the interpolation in $L^p$ spaces, we derive
\begin{equation}\label{5.9.3}v_n\rightarrow v\ \text{in}\ L^p(\mathbb{R}^N), \ \text{for any } p\in[2,\infty).\end{equation}
Note that $v_n$ solves
$$(-\Delta)^\alpha v_n+v_n=\chi_n,\ \ \ \ x\in \mathbb{R}^N,$$
where $\chi_n=f(v_n)-V(\epsilon_nx+\epsilon_ny_n)v_n+v_n.$ Denote
$\chi=f(v)-V(y_0)v+v.$ By (\ref{5.9.3}) and (\ref{2.2}) we know
$\chi_n\rightarrow \chi$ {in} $L^p(\mathbb{R}^N)$, for any $p\in[2,\infty)$.
Then there is $\kappa>0$ such that
$|\chi_n|_\infty\leq \kappa$, {for all} $n\in\mathbb{N}.$
Taking into account the results in \cite{FQT}, we know
$$v_n(x)=(\mathcal{K}*{\chi_n})(x)=\int_{\mathbb{R}^N}\mathcal{K}(x-y)\chi_n(y)dy,$$
where $\mathcal{K}$ is the  Bessel kernel. Then arguing as in the proof of \cite[Lemma 2.6]{AM}, (\ref{4.18.0}) holds true. \ \ \ $\Box$

Now we estimate the decay properties of $v_n$.
\begin{lemma}\label{l5.9.1}
There exists $C>0$ such that
$${v}_{n}(x)\leq C (1+|x|^{N+2\alpha})^{-1},\ \ \forall x\in\mathbb{R}^N,$$
uniformly in $n\in\mathbb{N}$.
\end{lemma}
{\bf Proof}: By \cite[Lemma 4.3]{FQT}, we know there exists $w$ such that
\begin{equation}\label{5.9.7}0<w(x)\leq{C}{(1+|x|^{N+2\alpha})^{-1}}\ \text{in}\ \mathbb{R}^N, \ \text{and}\ (-\Delta)^\alpha w +\frac{V_0}{2}w\geq0 \ \text{in}\ \mathbb{R}^N\backslash{B_{R_1}(0)},\end{equation}
for some $R_1>0$. In view of (\ref{4.18.0}), there exists $R_2>0$ such that $f(v_n)\leq\frac{V_0}{2}v_n$ for all $|x|\geq R_2$. Hence,
\begin{equation}\label{5.9.9}(-\Delta)^\alpha v_n +\frac{V_0}{2}v_n\leq f(v_n)-\frac{\tilde{V}_{\epsilon_n}(x)}{2}v_n\leq0, \ \text{in}\ \mathbb{R}^N\backslash{B_{R_2}(0)}.\end{equation}
Choose $R_3=\max\{R_1,R_2\}$, and we set
\begin{equation}\label{5.9.10} a=\inf_{B_{R_3}(0)}w>0, \ \text{and}\ Z_n:=(b+1)w-av_n,\end{equation}
where $b=\sup_{n\in\mathbb{N}}|v_n|_\infty<\infty$. Now we show
$Z_n\geq 0 $ {in} $\mathbb{R}^N.$
Argue by contradiction suppose that there is $x^k_{n}\subset\mathbb{R}^N$ such that
\begin{equation}\label{5.9.12}\inf_{x\in\mathbb{R}^N}Z_n(x)
=\lim_{k\rightarrow\infty}Z_n(x^k_{n})< 0. \end{equation}
In view of  $\lim_{|x|\rightarrow\infty}w(x)=0$ and (\ref{4.18.0})
we have
$\lim_{|x|\rightarrow\infty}Z_n(x)=0$ uniformly in $n\in\mathbb{N}$. Consequently, $\{x^k_{n}\}$ is bounded, and we may assume $x^k_{n}\rightarrow x^*_{n}$ for some $x^*_{n}\in\mathbb{R}^N$. Hence,
(\ref{5.9.12}) becomes
\begin{equation}\label{5.9.13}\inf_{x\in\mathbb{R}^N}Z_n(x)=Z_n(x^*_{n})< 0. \end{equation}
From the minimality property of $x^*_{n}$ and the integral representation of the fractional Laplace
of $Z_n$ at the point $x^*_{n}$,  it follows that
\begin{equation}\label{5.9.14}(-\Delta)^\alpha Z_n(x^*_{n})=C(N,\alpha)\int_{\mathbb{R}^N}
\frac{2Z_n(x^*_{n})-Z_n(x^*_{n}+y)-Z_n(x^*_{n}-y)}{|y|^{N+2\alpha}}dy\leq0.\end{equation}
By (\ref{5.9.10}) we get
$$Z_n(x)\geq (b+1)w-ba>0 \ \ \text{in}\ B_{R_3}(0).$$
Then $x^*_{n}\not\in {B_{R_3}(0)}$. From
(\ref{5.9.7}) and (\ref{5.9.9}) we know
$(-\Delta)^\alpha Z_n +\frac{V_0}{2}Z_n\geq0$\ {in}\ $\mathbb{R}^N\backslash{B_{R_3}(0)}.$
However, by (\ref{5.9.13}) and (\ref{5.9.14}) we have
$(-\Delta)^\alpha Z_n(x^*_{n}) +\frac{V_0}{2}Z_n(x^*_{n})<0.$
This is a contradiction. So $Z_n(x)\geq0$ in $\mathbb{R}^N$.
Using (\ref{5.9.10}) and (\ref{5.9.7}), the conclusion follows immediately.
\ \ \ \ \ $\Box$

{\bf Proof of Theorem 1.1.} Going back to (\ref{1.1}) with the variable substitution: $x\rightarrow \frac{x}{\epsilon}$, Lemma \ref{l4.7} implies that the equation (\ref{1.1}) has at
least $cat_{\mathcal{V}_\delta}(\mathcal{V})$ positive solutions for small $\epsilon$. Set $\epsilon_n\rightarrow0$. Let $w_n$ be the positive solutions of (\ref{1.1}) obtained above. Then $u_n(x)=w_n(\frac{x}{\epsilon_n})$ and $v_n(x)=u_n(x+y_n)$ are the associated solutions of (\ref{2.1}) and (\ref{4.7}) respectively. If $b_n$ denotes a maximum point of $v_n$, then from Lemma \ref{l5.9} it follows that $b_n\in B_{R}(0)$ for some $R>0$. Thus, the global maximum point of $w_n$ is $x_n:=\epsilon_nb_n+\epsilon_ny_n$.
Notice that $\epsilon_ny_n\rightarrow y_0\in \mathcal{V}$ using Lemma \ref{l4.5}. Then $x_n\rightarrow y_0$.
In addition, Lemma \ref{l5.9.1} and $|b_n|\leq R$ imply
$$\aligned w_n(x)&={v}_n(\frac{x}{\epsilon_n}-\frac{x_n}{\epsilon_n}-b_n)\leq \frac{C\epsilon^{N+2\alpha}_n}{\epsilon^{N+2\alpha}_n+|{x-x_n}|^{N+2\alpha}}.\endaligned$$
This ends the proof.\ \ \ \ $\Box$

\section{Proof of Theorem 1.2}
In this section, we are in a position to prove Theorem 1.2, and we always assume that (V$_1$), (V$_2$) and (f$_1$)-(f$_5$) are satisfied.
\subsection{Multiplicity of positive solutions}

For $a>0$, let $C_a(a^j)$ denote the hypercube $\Pi^N_{i=1}(a^j_i-a, a^j_i+a)$ centered at $a^j=(a^j_1,a^j_2,...,a^j_N)$, $j=1,...,k$. Let $\overline{C_a(a^j)}$ and $\partial C_a(a^j)$ denote the closure and the boundary of $C_a(a^j)$ respectively. Set $C^j_a\equiv C_a(a^j)$. By condition (V$_2$), we can choose numbers $L,l>0$ such that $\overline{C_a(a^j)}$ are disjoint, and $V(x)>V(a^j)$ for all $x\in \partial C_a(a^j)$, $j=1,...,k$, and $\cup^k_{j=1}C_l(a^j)\subset\Pi^N_{i=1}(-L, L)$. Define $\Psi_\epsilon\in C(\mathbb{R}, \mathbb{R})$, $h_\epsilon\in C(H^\alpha(\mathbb{R}^N),\mathbb{R}^N)$ by
\begin{equation}\label{4.4.1}\aligned\Psi_\epsilon(t)=\left\{ \begin{array}{ll}
-\frac{2L}{\epsilon},& \hbox{ $t<-\frac{2L}{\epsilon}$;}\\
 t, & \hbox{ $-\frac{2L}{\epsilon}\leq t\leq\frac{2L}{\epsilon}$;}\\
\frac{2L}{\epsilon}, & \hbox{ $t>\frac{2L}{\epsilon}$,}
\end{array}
\right.\endaligned\end{equation}
\begin{equation}\label{4.4.2}h^i_\epsilon(u)=\frac{\int_{\mathbb{R}^N}\Psi_\epsilon(x_i)|u|^pdx}{\int_{\mathbb{R}^N}|u|^pdx},\ i=1,2,...,N,\quad p\in[2,2^*_\alpha],\end{equation}
$$h_\epsilon(u)=(h^1_\epsilon(u),h^2_\epsilon(u),...,h^N_\epsilon(u)).$$
 Let $C^j_{{l}/{\epsilon}}\equiv C_{{l}/{\epsilon}}(\frac{a^j}{\epsilon})$ and
$$M^j_\epsilon=\{u\in M_\epsilon
: u\geq0\ {\text{and}}\  h_\epsilon(u)\in C^j_{{l}/{\epsilon}}\},\quad \alpha^j_\epsilon=\inf_{u\in M^j_\epsilon}I_\epsilon(u),$$
$$\partial M^j_\epsilon=\{u\in M_\epsilon
: u\geq0\ {\text{and}}\ h_\epsilon(u)\in \partial C^j_{{l}/{\epsilon}}\},\quad \tilde{\alpha}^j_\epsilon=\inf_{u\in \partial M^j_\epsilon}I_\epsilon(u),$$
for $j=1,...,k$.
We shall verify that $M^j_\epsilon$ and $\partial M^j_\epsilon$ are non-empty for $j=1,...,k$.  We remark that if  there exists a unique $t>0$ such that $tu\in M_\epsilon$ for any $u\in H^\alpha(\mathbb{R}^N)\backslash\{0\}$, then it is easy to see that $M^j_\epsilon$ and $\partial M^j_\epsilon$ are non-empty. Indeed, one can construct a function $u_0$ satisfying $h_\epsilon(u_0)\in C^j_{{l}/{\epsilon}}$ or $\partial C^j_{{l}/{\epsilon}}$. Moreover, there exists $t_0>0$ such that $t_0u_0\in M_\epsilon$. Then $t_0u_0\in C^j_{{l}/{\epsilon}}$ or $\partial C^j_{{l}/{\epsilon}}$.
However, in the current paper, for any $u\in H^\alpha(\mathbb{R}^N)\setminus\{0\}$, there may not exist $t>0$ such that $tu\in M_\epsilon$. Hence, it seems that it is not easy to verify that $M^j_\epsilon$ and $\partial M^j_\epsilon$ are non-empty.

\begin{lemma}\label{l5.1} There exists $\epsilon_0>0$ such that for any $\epsilon\in (0,\epsilon_0)$, $M^j_\epsilon$ and $\partial M^j_\epsilon$ are non-empty for all $j=1,...,k$.\end{lemma}
{\bf Proof}: Choose an even and nonnegative function $u\in C^\infty_0(\mathbb{R}^N)$ such that $|u|_2=1$. Then there exists $u_t(x):=t^{\frac{N}{2}}u(tx)$ with $t>0$ such that $|u_t|_2=1$ and $[ u_t]^2_\alpha\rightarrow0$ as $t\rightarrow0$.
Using $|V|_\infty<l_0$ in (f$_3$), there exists an  even and nonnegative function $w:=u_{t_0}\in C^\infty_0(\mathbb{R}^N)$ with small $t_0>0$  such that
\begin{equation}\label{6.2.3}[ w]^2_\alpha+|V|_\infty|w|^2_2-l_0|w|^2_2<0.\end{equation}
For small $\epsilon>0$ satisfying $2\sqrt{\epsilon}<1$, define $\chi_\epsilon\in C^1(\mathbb{R}^N,\mathbb{R}^+)$ by
$\chi_\epsilon(x)=
1$ if $|x|<\frac{1}{2\sqrt{\epsilon}}-1$, $\chi_\epsilon(x)=
0$ if $|x|>\frac{1}{2\sqrt{\epsilon}}$,
and $|\nabla \chi_\epsilon|\leq 2$ in $\mathbb{R}^N$. Take
$v^j_\epsilon(x)=w\bigl(x-\frac{a^j}{\epsilon}\bigr)\chi_\epsilon\bigl(x-\frac{a^j}{\epsilon}\bigr),$
for all $x\in\mathbb{R}^N$ and $j=1,...,k$.
By  Lebesgue's convergence theorem, we have
\begin{equation*} \aligned &\lim_{\epsilon\rightarrow0}[ v^j_\epsilon]^2_\alpha=[w]^2_\alpha,\quad \lim_{\epsilon\rightarrow0}|v^j_\epsilon|^2_2=|w|^2_2.\endaligned\end{equation*}
Then \begin{equation}\label{7.2}[ v^j_\epsilon]^2_\alpha+\int_{\mathbb{R}^N}V(\epsilon x)(v^j_\epsilon)^2dx-l_0|v^j_\epsilon|^2_2\leq[ w]^2_\alpha+|V|_\infty|w|^2_2-l_0|w|^2_2+o_\epsilon(1).\end{equation}
By (\ref{6.2.3}),  there exists $\epsilon_1>0$ such that $\epsilon\in (0,\epsilon_1)$, $v^j_\epsilon\in \Theta_\epsilon$.
Below we show $h_\epsilon(v^j_\epsilon)\in C^j_{{l}/{\epsilon}}$. Indeed,
from the definition of $\chi_\epsilon$, (\ref{4.4.1}) and (\ref{4.4.2}) it follows that
\begin{equation}\label{6.2.1}\aligned h^i_\epsilon( v^j_\epsilon)=&\frac{\int_{C^j_{1/{2\sqrt{\epsilon}}}(0)}\Psi_\epsilon(x_i+\frac{a^j}{\epsilon})
|w|^p\chi^p_\epsilon dx}
{\int_{C^j_{1/{2\sqrt{\epsilon}}}(0)}
|w|^p\chi^p_\epsilon dx}
=\frac{\int_{C^j_{1/{2\sqrt{\epsilon}}}(0)}(\frac{a^j_i}{\epsilon}+x_i)
|w|^p\chi^p_\epsilon dx}
{\int_{C^j_{1/{2\sqrt{\epsilon}}}(0)}|w|^p\chi^p_\epsilon dx}=\frac{a^j_i}{\epsilon},\endaligned\end{equation}
where we have restricted $\frac1{2\sqrt{\epsilon}}<\frac{l}{\epsilon}$ and used the fact that  $w$ is even.
Therefore, there exists $\epsilon_2<\min\{\epsilon_1,4l^2\}$ such that for $\epsilon\in (0,\epsilon_2)$, $h_\epsilon( v^j_\epsilon)\in C^j_{{l}/{\epsilon}}$, {and} $v^j_\epsilon\in \Theta_\epsilon.$ In view of Lemma \ref{l2.2} (ii), there exists $t^j_\epsilon>0$ such that $t^j_\epsilon v^j_\epsilon\in M_\epsilon
$.  So $t^j_\epsilon v^j_\epsilon\in M^j_\epsilon$.

Now we use similar argument to show  $\partial M^j_\epsilon\neq\emptyset$.
Denote $(v^j_\epsilon)^*(x)=v^j_\epsilon(x-\frac{l}{\epsilon}).$ As (\ref{7.2}) we infer
$$[ (v^j_\epsilon)^*]^2_\alpha+\int_{\mathbb{R}^N}V(\epsilon x)((v^j_\epsilon)^*)^2dx-l_0|(v^j_\epsilon)^*|^2_2\leq[ w]^2_\alpha+|V|_\infty|w|^2_2-l_0|w|^2_2+o_\epsilon(1).$$
Then there exists $\epsilon_3>0$ such that $\epsilon\in (0,\epsilon_3)$, $(v^j_\epsilon)^*\in \Theta_\epsilon$. Moreover, as (\ref{6.2.1}) we deduce
\begin{equation*}\aligned h^i_\epsilon( (v^j_\epsilon)^*)
=&\frac{\int_{C^j_{1/{2\sqrt{\epsilon}}}(0)}
(x_i+\frac{a^j_i}{\epsilon}+\frac{l}{\epsilon})
|w(x)|^p\chi^p_\epsilon(x)dx}
{\int_{C^j_{1/{2\sqrt{\epsilon}}}(0)}|w(x)|^p\chi^p_\epsilon(x)dx}=\frac{a^j_i}{\epsilon}+\frac{l}{\epsilon},
\endaligned\end{equation*}
where we have restricted $\frac1{2\sqrt{\epsilon}}<\frac{l}{\epsilon}$, used the fact that $2l\leq L$ and $w$ is even.
Then
for any $\epsilon<\min\{\epsilon_3,4l^2\}$, $h_\epsilon( (v^j_\epsilon)^*)\in \partial C^j_{{l}/{\epsilon}}$ {and}\ $( v^j_\epsilon)^*\in \Theta_\epsilon.$ According to Lemma \ref{l2.2} (ii),  there exists $(t^j_\epsilon)^*>0$ such that $(t^j_\epsilon)^*( v^j_\epsilon)^*\in M_\epsilon$ and so $(t^j_\epsilon)^*( v^j_\epsilon)^*\in \partial M^j_\epsilon$.
\ \ \ \ $\Box$
\begin{Remark}For simplicity, we select the same test function $w$ in the proof of $M^j_\epsilon\neq\emptyset$ as in that of $\partial M^j_\epsilon\neq\emptyset$. Actually, we can choose another test function $w_0$ as in the proof of the following lemma to show $M^j_\epsilon\neq\emptyset$. \end{Remark}
\begin{lemma}\label{l5.2} $\lim_{\epsilon\rightarrow0}c_\epsilon=\lim_{\epsilon\rightarrow0}\alpha^j_\epsilon= c_{V_0},$
for all $j=1,...,k$. \end{lemma}
{\bf Proof}:  Take $w_0$ be a positive ground state of the equation (\ref{1.0.0}).
Define $\tilde{v}^j_\epsilon(x)=w_0\bigl(x-\frac{a^j}{\epsilon}\bigr)\chi_\epsilon\bigl(x-\frac{a^j}{\epsilon}\bigr).$
for all $x\in\mathbb{R}^N$ and $j=1,...,k$, where $\chi_\epsilon$ is given as in Lemma \ref{l5.1}.
 Moreover, Lebesgue's convergence theorem implies
\begin{equation}\label{6.7.1}\aligned &\lim_{\epsilon\rightarrow0}\int_{\mathbb{R}^N}V(\epsilon x)(\tilde{v}^j_\epsilon)^2dx=V_0|w_0|^2_2, \quad \lim_{\epsilon\rightarrow0}[ \tilde{v}^j_\epsilon]^2_\alpha=[w_0]^2_\alpha,\\& \lim_{\epsilon\rightarrow0}|\tilde{v}^j_\epsilon|^2_2=|w_0|^2_2,\ \lim_{\epsilon\rightarrow0}\int_{\mathbb{R}^N}f(\tilde{v}^j_\epsilon)\tilde{v}^j_\epsilon dx=\int_{\mathbb{R}^N}f(w_0)w_0dx.\endaligned\end{equation}
Similar to (\ref{4.7.1}), we have $\tilde{v}^j_\epsilon\in \Theta_\epsilon$ for small $\epsilon$. By restricting $\frac1{2\sqrt{\epsilon}}<\frac{l}{\epsilon}$, as (\ref{6.2.1}) we get
\begin{equation*} h^i_\epsilon( \tilde{v}^j_\epsilon)=\frac{\int_{C^j_{1/{2\sqrt{\epsilon}}}(0)}(\frac{a^j_i}{\epsilon}+x_i)
|w_0(x)|^p\chi^p_\epsilon(x)dx}
{\int_{C^j_{1/{2\sqrt{\epsilon}}}(0)}|w_0(x)|^p\chi^p_\epsilon(x)dx}\in\Bigl(\frac{a^j_i}{\epsilon}-\frac1{2\sqrt{\epsilon}}, \frac{a^j_i}{\epsilon}+\frac1{2\sqrt{\epsilon}}\Bigr)\subset C^j_{\frac{l}{\epsilon}},\end{equation*}
for small enough $\epsilon$. According to Lemma \ref{l2.2} (ii) there exists $t^j_\epsilon>0$ such that $t^j_\epsilon \tilde{v}^j_\epsilon\in M_\epsilon$ and then $t^j_\epsilon \tilde{v}^j_\epsilon\in M^j_\epsilon$ for  $j=1,...,k$.
 Similar to (\ref{4.9.0}), we deduce that  $t^j_\epsilon\rightarrow1$ as $\epsilon\rightarrow0$.
From (\ref{6.7.1}) it follows that
$I_\epsilon(t^j_\epsilon \tilde{v}^j_\epsilon)+o_\epsilon(1)=I_{V_0}(w_0)=c_{V_0}.$
Observe that $I_\epsilon(t^j_\epsilon \tilde{v}^j_\epsilon)\geq \alpha^j_\epsilon$. Then $\limsup_{\epsilon\rightarrow0}\alpha^j_\epsilon\leq c_{V_0}$. Using (\ref{2.5}) and (\ref{3.2}) we know $\liminf_{\epsilon\rightarrow0}c_\epsilon\geq c_{V_0}$. In addition, it is clear that
$c_\epsilon\leq \alpha^j_\epsilon$.  Therefore
$\lim_{\epsilon\rightarrow0}c_\epsilon=\lim_{\epsilon\rightarrow0}\alpha^j_\epsilon= c_{V_0},$
for all $j=1,...,k$. \ \ \ \ $\Box$

\begin{lemma}\label{l5.3} $\liminf_{\epsilon\rightarrow0}\tilde{\alpha}^j_\epsilon> c_{V_0}$ for all $j=1,...,k$.\end{lemma}
{\bf Proof}: We argue by contradiction. Since $\tilde{\alpha}^j_\epsilon\geq c_\epsilon\geq c_{V_0}$, we may assume that there exists $\epsilon_n\rightarrow0$ such that $\tilde{\alpha}^j_{\epsilon_n}\rightarrow c_{V_0}$.
By the definition of $\tilde{\alpha}^j_{\epsilon_n}$, there exists $u_n\in M_{\epsilon_n}$ satisfying
\begin{equation}\label{6.2.2}h_{\epsilon_n}(u_n)\in \partial C^j_{\frac{l}{\epsilon_n}}, \quad I_{\epsilon_n}(u_n)\rightarrow c_{V_0}.\end{equation}
Since $u_n\in M_{\epsilon_n}$, we have $u_n\in \Theta_{\epsilon_n}\subset \Theta_{V_0}$. Lemma \ref{l3.1} (ii) implies that for some $t_n>0$, $t_nu_n\in  M_{V_0}$. Then
$$c_{V_0}+o_n(1)= I_{\epsilon_n}(u_n)\geq I_{\epsilon_n}(t_nu_n)\geq I_{V_0}(t_nu_n)\geq c_{V_0}.$$
Therefore,
$t_n\rightarrow1$, {and} \ $\int_{\mathbb{R}^N}V(\epsilon_n x)u^2_ndx=V_0|u_n|^2_2+o_n(1).$
Then $\bar{u}_n:=t_nu_n\in M_{V_0}$ is a minimizing sequence of $I_{V_0}$.
By Ekeland variational principle, we suppose
$(I_{V_0}|_{M_{V_0}})'(\bar{u}_n)\rightarrow0$. From Lemma \ref{l3.1} (iii) we infer that
$I'_{V_0}(\bar{u}_n)\rightarrow0$, $\{\bar{u}_n\}$ is bounded and  nonvanishing. Then there exists $y_n\in\mathbb{R}^N$ such that $v_n(\cdot):=\bar{u}_n(\cdot+y_n)\rightharpoonup v\neq0$ in $H^\alpha(\mathbb{R}^N)$.  Since $\{v_n\}\subset M_{V_0}$ is also a (PS)$_{c_{V_0}}$ sequence of $I_{V_0}$, from Lemma \ref{l4.4} we know
 $\|v_n-v\|\rightarrow0$, i.e. $\|\bar{u}_n(\cdot+y_n)-v(\cdot)\|\rightarrow0$. In view of $t_n\rightarrow1$,
we have  $\|{u_n}(\cdot+y_n)-v(\cdot)\|\rightarrow0$ and then $|{u_n}(\cdot+y_n)-v(\cdot)|_p\rightarrow0$ for any $p\in(2,2^*_\alpha)$. Therefore, for any $\delta>0$ there exists $R_\delta>0$ such that
\begin{equation}\label{7.6}\frac{\int_{|x|\geq R_\delta}|u_n(x+y_n)|^pdx}{\int_{\mathbb{R}^N}|u_n(x+y_n)|^pdx}\leq\delta.\end{equation}
By (V$_2$) we may suppose that, there exists $\nu>0$ such that
\begin{equation}\label{7.7}V(x)>V_0 \quad \text{for }\ x\in \bar{C}^j_{l+\nu}\backslash{{C}^j_{l-\nu}}.\end{equation}
 Below we consider the sequence $\{\epsilon_n y_n\}$. Passing to a subsequence, we may assume that one of the following cases occurs.\\
(1) $\epsilon_n y_n\rightarrow y_0\in\bar{C}^j_{l+\nu}\backslash{{C}^j_{l-\nu}}$;\\
(2) $\{\epsilon_ny_n\}\subset\bar{C}^j_{l-\nu}$;\\
(3) $\{\epsilon_ny_n\}\subset\mathbb{R}^N\backslash{{C}^j_{l+\nu}}$ and $\{\epsilon_n y_n\}$ is bounded;\\
(4) $\epsilon_n(y_n)_i\rightarrow\infty$ as $n\rightarrow\infty$ for some $i\in{1,2,..., N}$.

If case (1) occurs, by (\ref{7.7}) we have $V(y_0)>V_0$.
Since $\|{u_n}(\cdot+y_n)-v\|\rightarrow0$ and (\ref{6.2.2}), we obtain
$$\aligned c_{V_0}= \liminf_{n\rightarrow\infty}I_{\epsilon_n}(u_n)=&\liminf_{n\rightarrow\infty}
\Bigl\{\frac12[ u_n(x+y_n)]^2_\alpha+\frac12\int_{\mathbb{R}^N}V(\epsilon_nx+\epsilon_ny_n)
u^2_n(x+y_n)dx\\&-
\int_{\mathbb{R}^N}F(u_n(x+y_n))dx\Bigr\}\\ \geq&\frac12[v]^2_\alpha+\frac {V(y_0)}{2}|v|^2_2-
\int_{\mathbb{R}^N}F(v)dx>I_{V_0}(v)\geq c_{V_0}.\endaligned$$
This is a contradiction.

 Now we consider case (2). From (\ref{4.4.2}) it follows that
\begin{equation}\label{7.8}\aligned &h^i_{\epsilon_n}(u_n)=&\frac{\bigl(\int_{|x|\leq R_\delta}+\int_{|x|\geq R_\delta}\bigr)\Psi_{\epsilon_n}((x+y_n)_i)|u_n(x+y_n)|^pdx}
{\int_{\mathbb{R}^N}|u_n(x+y_n)|^pdx}
,\endaligned\end{equation}
where $i=1,2,...,N$. For $|x_i|\leq R_\delta$ we have
 $$(x+y_n)_i\in \Bigl(\frac{a^j_i-(l-\nu)}{\epsilon_n}-R_\delta,\frac{a^j_i+(l-\nu)}{\epsilon_n}+R_\delta \Bigr)\subset\Bigl(-\frac{2L}{\epsilon_n},\frac{2L}{\epsilon_n}\Bigr),$$
 for large $n$. According to (\ref{7.6}) and the definition of $\Psi_{\epsilon_n}$ we infer
$$\Bigl(\frac{a^j_i-(l-\nu)}{\epsilon_n}-R_\delta\Bigr)(1-\delta)-\frac{2L}{\epsilon_n}\cdot\delta\leq h^i_{\epsilon_n}(u_n)
\leq\Bigl(\frac{a^j_i+(l-\nu)}{\epsilon_n}+R_\delta\Bigr)+\frac{2L}{\epsilon_n}\cdot\delta.$$
Choose $\delta>0$ and $\nu>\delta$ small enough such that
$$\frac{a^j_i-l}{\epsilon_n}\leq h^i_{\epsilon_n}(u_n)\leq \frac{a^j_i+l}{\epsilon_n}\quad \text{for sufficiently large }\ n,$$
contradicting with $h_{\epsilon_n}(u_n)\in \partial C^j_{{l}/{\epsilon_n}}$.

In case (3), we may assume $\epsilon_ny_n\rightarrow \hat{y}$. Then $\hat{y}_i\geq a^j_i+l+\nu$ for some $i$ and $\epsilon_n(y_n)_i\geq a^j_i+l+\frac{\nu}{2}$ for large $n$. For $|x_i|\leq R_\delta$, we then have
$x_i+(y_n)_i\geq\frac{a^j_i+l+\frac{\nu}{2}}{\epsilon_n}-R_\delta.$
So from (\ref{7.8}) we get
$$h^i_{\epsilon_n}(u_n)\geq\Bigl(\frac{a^j_i+l+\frac{\nu}{2}}{\epsilon_n}-R_\delta\Bigr)(1-\delta)-\frac{2L}{\epsilon_n}\cdot\delta
>\frac{a^j_i+l}{\epsilon_n},$$
for sufficiently small $\delta>0$, $\nu>\delta$ and $n$ large enough. This is impossible since $h_{\epsilon_n}(u_n)\in \partial C^j_{{l}/{\epsilon_n}}$.

Case (4) is excluded by assuming $\epsilon_n(y_n)_i>2L$ or $\epsilon_n(y_n)_i<-2L$ for some $i$ and  large $n$, and using a similar argument to that of Case (3). \ \ \ \ $\Box$

In view of Lemmas \ref{l5.1}, \ref{l5.2} and \ref{l5.3} and $c_{V_0}<c_{V_\infty}$, we can find a small positive number $\xi$ and take $\epsilon^*<\epsilon_0$ such that, for every $\epsilon\in (0,\epsilon^*)$ and $j=1,...,k$,
\begin{equation}\label{7.10} c_\epsilon\leq \alpha^j_\epsilon< c_{V_0}+\xi<\min\{2c_{V_0},\tilde{\alpha}^j_\epsilon, c_{V_\infty}\}. \end{equation}

In order to look for PS sequences, inspired by \cite{CN} we obtain the following result:
\begin{lemma}\label{l5.4} For every $u\in M^j_\epsilon$, $j=1,...,k$ and $\epsilon\in (0,\epsilon^*)$, there exist $\delta>0$ and a differentiable function $t=t(w)>0$, where $w\in H^\alpha(\mathbb{R}^N)$ with $\|w\|_\epsilon\leq\delta$, such that $t(0)=1$, the function $z=t(w)(u-w)\in M^j_\epsilon$ for all $w\in H^\alpha(\mathbb{R}^N)$ with $\|w\|_\epsilon\leq\delta$ and
\begin{equation}\label{5.12.0}\langle t'(0),\phi\rangle=\frac{2\langle u,\phi\rangle_\epsilon-\int_{\mathbb{R}^N}f'(u)u\phi dx-\int_{\mathbb{R}^N}f(u)\phi dx}{\|u\|_\epsilon-\int_{\mathbb{R}^N}f'(u)u^2 dx},\end{equation}
for any $\phi\in H^\alpha(\mathbb{R}^N)$.
\end{lemma}
{\bf Proof}: For $u\in M^j_\epsilon$, define a function $H_u:\mathbb{R}\times H^\alpha(\mathbb{R}^N)\rightarrow\mathbb{R}$ by
$$H_u(t,w)=\langle I'_\epsilon(t(u-w)),u-w\rangle=t\|u-w\|^2_\epsilon-\int_{\mathbb{R}^N}f(t(u-w))(u-w)dx.$$
Then $H_u(1,0)=0$ and using Lemma \ref{l2.1} (iii) there holds
$$\frac{\partial H_u}{\partial t}\Big|_{(1,0)}=\|u\|^2_\epsilon-\int_{\mathbb{R}^N}f'(u)u^2dx=\int_{\mathbb{R}^N}(f(u)u-f'(u)u^2)dx<0.$$
According to the implicit function theorem, there exist $\delta>0$ and a differentiable function $t(w)$
defined for $w\in H^\alpha(\mathbb{R}^N)$ with $\|w\|_\epsilon<\delta$ such that $t(0)=1$, the function $t$ also satisfies (\ref{5.12.0}),
and $H_u(t(w),w)=0$ for all $w\in H^\alpha(\mathbb{R}^N)$ with $\|w\|_\epsilon<\delta$. Then $t(w)(u-w)\in M_\epsilon
$.
Moreover, the continuity of the maps $h_\epsilon$ and $t$ imply that
$h_\epsilon(t(w)(u-w))\in C^j_{{l}/{\epsilon}},$
for small $\delta$. Therefore, $t(w)(u-w)\in M^j_\epsilon$  for all $w\in H^\alpha(\mathbb{R}^N)$ with $\|w\|_\epsilon<\delta$. \ \ \ \ $\Box$

\begin{lemma}\label{l5.5} For all $\epsilon\in (0,\epsilon^*)$ and $j=1,...,k$, there exists a minimizing sequence $u^j_n\in M^j_\epsilon$ such that
$I_\epsilon(u^j_n)\rightarrow\alpha^j_\epsilon$ and $I'_\epsilon(u^j_n)\rightarrow 0$ {in} $(H^{\alpha}(\mathbb{R}^N))^*$
as $n\rightarrow\infty$.\end{lemma}
{\bf Proof}: Denote $\bar{M}^j_\epsilon={M}^j_\epsilon\cup \partial{M}^j_\epsilon$.
By (\ref{7.10}) we get $\alpha^j_\epsilon=\inf\{I_\epsilon(u)|u\in \bar{M}^j_\epsilon\}$. Applying the Ekeland variational
principle, there exists a minimizing sequence $\{u_n\}\subset {M}^j_\epsilon$ such that
\begin{equation}\label{7.11}I_\epsilon(u_n)<\alpha^j_\epsilon+\frac1n,\ I_\epsilon(v)\geq I_\epsilon(u_n)-\frac1n\|u_n-v\|_\epsilon,\end{equation}
for all $v\in \bar{{M}}^j_\epsilon$. Applying Lemma \ref{l5.4} with $u=u_n$, there exists  $t^*_n(w)$ defined for $w\in H^\alpha(\mathbb{R}^N)$ and $\|w\|_\epsilon<\delta_n$ for some $\delta_n>0$
such that $t^*_n(w)(u_n-w)\in {M}^j_\epsilon$. Choose $0<\sigma<\delta_n$. Let $u\in H^\alpha(\mathbb{R}^N)$, $u\not\equiv0$ and $w_\sigma=\frac{\sigma u}{\|u\|_\epsilon}$. Fix $n$ and set $z_{\sigma}:=t^*_n(w_\sigma)(u_n-w_\sigma)$. Observe that $z_\sigma\in M^j_\epsilon$.  From (\ref{7.11}) and mean value theorem it follows that
$$-\frac1n\|z_\sigma-u_n\|_\epsilon \leq I_\epsilon(z_\sigma)-I_\epsilon(u_n)=\langle I'_\epsilon(u_n),z_\sigma-u_n\rangle+o(\|z_\sigma-u_n\|_\epsilon).$$
Then
\begin{equation}\label{7.12}\langle I'_\epsilon(u_n),-w_\sigma\rangle+(t^*_n(w_\sigma)-1)\langle I'_\epsilon(u_n),u_n-w_\sigma\rangle\geq-\frac1n\|z_\sigma-u_n\|_\epsilon+o(\|z_\sigma-u_n\|_\epsilon).\end{equation}
Using $z_\sigma\in {M}^j_\epsilon$ again, we have $\langle I'_\epsilon(z_\sigma),u_n-w_\sigma\rangle=0.$
Then (\ref{7.12}) implies
\begin{equation}\label{7.13}\langle I'_\epsilon(u_n),\frac{u}{\|u\|_\epsilon}\rangle\leq \frac{\|z_\sigma-u_n\|_\epsilon}{n\sigma}+\frac{o(\|z_\sigma-u_n\|_\epsilon)}{\sigma}
+\frac{(t^*_n(w_\sigma)-1)}{\sigma}\langle I'_\epsilon(u_n)-I'_\epsilon(z_\sigma),u_n-w_\sigma\rangle.\end{equation}
From the definition of $z_\sigma$, we get
\begin{equation}\label{5.2.0}\|z_\sigma-u_n\|_\epsilon\leq\sigma+C(|t^*_n(w_\sigma)-1|).\end{equation}
According to Lemma \ref{l5.4}, there holds
$$\langle (t^*_n)'(0),\phi\rangle=\frac{2\langle u_n,\phi\rangle_\epsilon-\int_{\mathbb{R}^N}f'(u_n)u_n\phi dx-\int_{\mathbb{R}^N}f(u_n)\phi dx}{\|u_n\|^2_\epsilon-\int_{\mathbb{R}^N}f'(u_n)u^2_n dx}.$$
In view of Lemma \ref{l5.2} and $I_\epsilon(u_n)\rightarrow\alpha^j_\epsilon$, as Lemma \ref{l2.3} we conclude that  $\{u_n\}$ is bounded in $H^\alpha(\mathbb{R}^N)$ and nonvanishing.  Then there exists $y_n\in\mathbb{R}^N$ such that
$u_n(\cdot+y_n)\rightharpoonup u_0\neq0$ in $H^\alpha(\mathbb{R}^N)$. As the argument in Lemma \ref{l2.4} we have $$\Bigl|2\langle u_n,\phi\rangle_\epsilon-\int_{\mathbb{R}^N}f'(u_n)u_n\phi dx-\int_{\mathbb{R}^N}f(u_n)\phi dx\Bigr|\leq C\|\phi\|,$$
and
$$\Bigl|\|u_n\|^2_\epsilon-\int_{\mathbb{R}^N}f'(u_n)u^2_n dx\Bigr|\geq\int_{\mathbb{R}^N}[f'(u_0)u^2_0 -f(u_0)u_0] dx>0.$$
Then $\|(t^*_n)'(0)\|_{(H^{\alpha}(\mathbb{R}^N))^*}\leq C$.
Hence,
$\lim_{\sigma\rightarrow0}\frac{|t^*_n(w_\sigma)-1|}{\sigma}
\leq C,$
and  $\lim_{\sigma\rightarrow0}\frac{\|z_\sigma-u_n\|_\epsilon}{\sigma}\leq C$ using (\ref{5.2.0}).
For fixed $n$, from (\ref{7.13}) we deduce
$\lim_{\sigma\rightarrow0}\langle I'_\epsilon(u_n),\frac{u}{\|u\|_\epsilon}\rangle\leq\frac{C}{n}.$
Letting $n\rightarrow\infty$, we know
$I'_\epsilon(u_n)\rightarrow0$ in $(H^{\alpha}(\mathbb{R}^N))^*$.\ \ \ \ $\Box$

\begin{lemma}\label{l5.7} Equation (\ref{2.1})
has at least $k$ positive solutions for each $\epsilon\in (0,\epsilon^*)$.\end{lemma}
{\bf Proof}: In view of Lemma \ref{l5.5}, for fixed $j=1,...,k$, there exists $\{u^j_n\}\subset M^j_\epsilon$ satisfying $$I_\epsilon(u^j_n)\rightarrow\alpha^j_\epsilon,\quad \text{and}\ I'_\epsilon(u^j_n)\rightarrow0.$$For fixed $j=1,...,k$ and $\epsilon\in (0,\epsilon^*)$,
 by (\ref{7.10}) and Lemma \ref{l4.3}, we may assume that $u^j_n\rightarrow u^j_0$ in $H^\alpha(\mathbb{R}^N)$. So $I'_\epsilon(u^j_0)=0$,  $I_\epsilon(u^j_0)=\alpha^j_\epsilon$ and $u^j_0\in M^j_\epsilon\cup\partial M^j_\epsilon$. Using $\alpha^j_\epsilon<\tilde{\alpha}^j_\epsilon$ and $I_\epsilon(u^j_0)=\alpha^j_\epsilon$, we have that $u^j_0\in M^j_\epsilon$.
It is clear that $u^j_0$ is non-negative, and  the
maximum principle implies $u^j_0$ is positive for $j=1,...,k$. Moreover,
$u^j_0$ ($j=1,...,k$) are different since $C^j_{{l}/{\epsilon}}$ are disjoint. Therefore, the equation (\ref{2.1})
has at least $k$ positive solutions.\ \ \ \ $\Box$

\subsection{Concentration of positive solutions}

\begin{lemma}\label{l5.8} Assume that the potential $V$ is uniformly continuous. Let $u^j_{\epsilon}$, $j=1,...,k$, be the positive solutions obtained in Lemma \ref{l5.7}. Then there is $y^j_{\epsilon}\in\mathbb{R}^N$
such that $\epsilon y^j_\epsilon\rightarrow a^j$ as
$\epsilon\rightarrow 0$. Setting $v^j_\epsilon(x):=u^j_\epsilon( x+y^j_\epsilon)$, then $v^j_\epsilon$ converges in $H^\alpha(\mathbb{R}^N)$ to a positive solution $v^j$ of the equation (\ref{1.0.0}).\end{lemma}
{\bf Proof}: Fix $j\in\{1,...,k\}$. Let $u^j_n$ be the positive solutions of problem (\ref{2.1}) with parameter $\epsilon_n\rightarrow0$ satisfying
$$u^j_n\in M^j_{\epsilon_n},\ I_{\epsilon_n}(u^j_n)=\alpha^j_{\epsilon_n}, \ \text{and}\ I'_{\epsilon_n}(u^j_n)=0.$$
From Lemma \ref{l5.2} we know $\lim_{n\rightarrow\infty}\alpha^j_{\epsilon_n}=c_{V_{0}}$. Similar to Lemma \ref{l2.3}, we infer that $\{u^j_n\}$ is bounded in $H^\alpha(\mathbb{R}^N)$ and non-vanishing. Then there exist $y^j_n\in \mathbb{R}^N$ and $\delta>0$ such that
\begin{equation}\label{5.4}\int_{B_1(y^j_n)}(u^j_n(x))^2dx\geq\delta.\end{equation}
Setting $v^j_n(x)=u^j_n(x+y^j_n)$ and $\tilde{V}_{\epsilon_n}(x)=V(\epsilon_n(x+y^j_n))$, we see that $v^j_n$ solves
$$(-\Delta)^\alpha u+\tilde{V}_{\epsilon_n}(x)u=f(u),\ \ \ \text{in}\ \mathbb{R}^N,$$
with energy functional
$$\tilde{I}_{\epsilon_n}(v^j_n)=\frac12[ {v^j_n}]^2_\alpha+\frac12\int_{\mathbb{R}^N}\tilde{V}_{\epsilon_n}(x)(v^j_n)^2dx-\int_{\mathbb{R}^N}F(v^j_n)dx.$$
Since $\{v^j_n\}$ is bounded in $H^\alpha(\mathbb{R}^N)$, from (\ref{5.4}) we assume $v^j_n\rightharpoonup v^j\neq0$ in $H^\alpha(\mathbb{R}^N)$.

 {\bf Claim 1}. The sequence $\{\epsilon_ny^j_n\}$ must be bounded.

  Otherwise if $\epsilon_ny^j_n\rightarrow\infty$, then we may suppose $V(\epsilon_ny^j_n)\rightarrow \nu_1 \geq V_\infty>V_{0}$. Since $V$ is uniformly continuous, for $R>0$ and $|x|\leq R$ we have
\begin{equation}\label{5.3}|\tilde{V}_{\epsilon_n}(x)-\nu_1|\leq|V(\epsilon_n(x+y^j_n))-V(\epsilon_ny^j_n)|+|V(\epsilon_ny^j_n)-\nu_1|\rightarrow0.\end{equation}
For each $\eta\in C^\infty_0(\mathbb{R}^N)$, we infer
\begin{equation*}\lim_{n\rightarrow\infty}\int_{\mathbb{R}^N}\tilde{V}_{\epsilon_n}(x)v^j_n\eta dx =\int_{\mathbb{R}^N}\nu_1 v^j\eta dx,\  \lim_{n\rightarrow\infty}\int_{\mathbb{R}^N}f(v^j_n)\eta dx =\int_{\mathbb{R}^N}f(v^j)\eta dx.
\end{equation*}
Then $v^j$ solves
\begin{equation}\label{7.2.2}
(-\Delta)^\alpha
v^j+\nu_1 v^j=f(v^j),\quad \text{in}\
\mathbb{R}^{N}.
\end{equation}
From Fatou's lemma, (\ref{2.3}) and Lemma \ref{l5.2} it follows that
\begin{equation}\label{4.3.5}
\begin{split}
c_{V_0}&=\liminf_{n\rightarrow\infty}\alpha^j_{\epsilon_{n}}
=\liminf_{n\rightarrow\infty}[I_{\epsilon_{n}}(u^j_{n})-\frac{1}{2}\langle I'_{\epsilon_{n}}(u^j_{n}),u^j_{n}\rangle]\\
&=\liminf_{n\rightarrow\infty}\bigl[\tilde{I}_{\epsilon_{n}}(v^j_{n})
-\frac{1}{2}\langle\tilde{I}'_{\epsilon_{n}}(v^j_{n}),v^j_{n}\rangle\bigr]
=\liminf_{n\rightarrow\infty}\int_{\mathbb{R}^{N}}
\bigl[\frac12f(v^j_{n})v^j_{n}-F(v^j_{n})\bigr]dx\\
&\geq\int_{\mathbb{R}^{N}}
\bigl[\frac12f(v^j)v^j-F(v^j)\bigr]dx={I}_{\nu_1}(v^j)-\frac12\langle{I}'_{\nu_1}(v^j),v^j\rangle\geq c_{\nu_1}.
\end{split}
\end{equation}
 However, by the fact that
$V_{0}<\nu_1$, Lemma \ref{l3.3} implies that $c_{V_{0}}< c_{\nu_1}$. This is a contradiction.
Hence $\{\epsilon_ny^j_n\}$ is bounded and suppose $\epsilon_ny^j_n\rightarrow y^j_0$. Then $v^j$ is a solution of the equation
(\ref{3.1}) with $a=V(y^j_0).$

{\bf Claim 2}.  $y^j_0\in \mathcal{V}=\{x\in \mathbb{R}^N: V(x)=V_0\}.$

If $y^j_0\not\in \mathcal{V}$, then by Lemma \ref{l3.3} we have $c_{V_{0}}<c_{V(y^j_0)}.$
However, repeating the argument of (\ref{4.3.5}), we have
$c_{V_{0}}\geq c_{V(y^j_0)}$. This is a contradiction. Then $y^j_0\in \mathcal{V}$ and so $V(y^j_0)=V_0$.

{\bf Claim 3}. $v^j_n$ converges strongly to $v^j$ in $H^\alpha(\mathbb{R}^N)$.

If $v^j_n\rightarrow v^j$ in $H^\alpha(\mathbb{R}^N)$, then we are done. Otherwise, set $w^j_{n,1}(x)=u^j_n(x)-v^j(x-y^j_n)$,
then $w^j_{n,1}(x+y^j_n)\rightharpoonup0$ and $w^j_{n,1}\not\rightarrow0$ in $H^\alpha(\mathbb{R}^N)$.
For any $\varphi\in H^\alpha(\mathbb{R}^N)$, note that $\nu_1=V(y^j_0)=V_0$,
by (\ref{5.3}) we obtain
$$\int_{\mathbb{R}^N}(\tilde{V}_{\epsilon_n}(x)-V_0)v^j
\varphi(x+y^j_n)dx=o_n(1)\|\varphi\|.$$
Similar to (\ref{4.2}) we have
\begin{equation}\label{5.5}\aligned\langle I'_{\epsilon_n}(w^j_{n,1}),\varphi\rangle&=\langle \tilde{I}'_{\epsilon_n}(w^j_{n,1}(x+y^j_n)),\varphi(x+y^j_n)\rangle\\&=\langle \tilde{I}'_{\epsilon_n}(v^j_{n}),\varphi(x+y^j_n)\rangle-\langle \tilde{I}'_{\epsilon_n}(v^j),\varphi(x+y^j_n)\rangle+o_n(1)\|\varphi\|\\
&=0-\langle I'_{V_0}(v^j),\varphi(x+y^j_n)\rangle+o_n(1)\|\varphi\|=o_n(1)\|\varphi\|.\endaligned\end{equation}
Then
\begin{equation}\label{5.1.1}I'_{\epsilon_n}(w^j_{n,1})\rightarrow0\  \text{in}\ (H^{\alpha}(\mathbb{R}^N))^*.\end{equation}
Using (\ref{5.3}) again we get
$\int_{\mathbb{R}^N}(\tilde{V}_{\epsilon_n}(x)-V_0)(v^j)^2dx\rightarrow0.$
As (\ref{4.3}) we conclude
$$\tilde{I}_{\epsilon_n}(v^j_n)=\tilde{I}_{\epsilon_n}(w^j_{n,1}(x+y^j_n))+\tilde{I}_{\epsilon_n}(v^j)+o_n(1)=\tilde{I}_{\epsilon_n}(w^j_{n,1}(x+y^j_n))+I_{V_0}(v^j)+o_n(1).$$
Then
\begin{equation}\label{5.1.2}I_{\epsilon_n}(u^j_n)-I_{\epsilon_n}(w^j_{n,1})-I_{V_0}(v^j)=o_n(1).\end{equation}
In view of (\ref{5.1.1}) and $w^j_{n,1}\not\rightarrow 0$ in $H^\alpha(\mathbb{R}^N)$ we get $\{w^j_{n,1}\}$ is non-vanishing. Hence, there exist $y^j_{n,2}\in\mathbb{R}^N$ and $\delta'>0$ such that
$\int_{B_1(y^j_{n,2})}(w^j_{n,1}(x))^2 dx>\delta'.$ Passing to a subsequence, we assume that $v^j_{n,1}(\cdot):=w^j_{n,1}(\cdot+y^j_{n,2})\rightharpoonup \tilde{v}^j_2\neq0$ in $H^{\alpha}(\mathbb{R}^N)$. By (\ref{5.1.1}) we have
 \begin{equation}\label{5.1.3}\tilde{\tilde{I}}'_{\epsilon_n}(v^j_{n,1})\rightarrow0\  \text{in}\ (H^{\alpha}(\mathbb{R}^N))^*, \end{equation}
with
$$\tilde{\tilde{I}}_{\epsilon_n}(v^j_{n,1})=\frac12[ {v^j_{n,1}}]^2_\alpha+\frac12\int_{\mathbb{R}^N}\tilde{\tilde{V}}_{\epsilon_n}(x)(v^j_{n,1})^2dx
-\int_{\mathbb{R}^N}F(v^j_{n,1})dx,$$ where $\tilde{\tilde{V}}_{\epsilon_n}(x)=V(\epsilon_n(x+y^j_{n,2}))$.  Denote
 $w^j_{n,2}(x)=w^j_{n,1}(x)-\tilde{v}^j_2(x-y^j_{n,2})$ and set $V(\epsilon_ny^j_{n,2})\rightarrow \nu_2$, where $\nu_2=V(z_0)$ if $\epsilon_ny^j_{n,2}$ is bounded and assume $\epsilon_ny^j_{n,2}\rightarrow z_0$, or else  $\nu_2\geq V_\infty$.
Since $V$ is uniformly continuous, we infer
$\tilde{\tilde{V}}_{\epsilon_n}(x)=V(\epsilon_nx+\epsilon_ny^j_{n,2})\rightarrow \nu_2,$
on any bounded set. Therefore, by (\ref{5.1.3}), similar to (\ref{7.2.2}) we get $\tilde{v}^j_2$ solves (\ref{3.1}) with $a=\nu_2$.
For any $\varphi\in H^\alpha(\mathbb{R}^N)$, arguing as (\ref{5.5}) we deduce
$$\aligned\langle {I}'_{\epsilon_n}(w^j_{n,2}),\varphi\rangle&=\langle \tilde{\tilde{I}}'_{\epsilon_n}(v^j_{n,1}),\varphi(x+y^j_{n,2})\rangle-\langle \tilde{\tilde{I}}'_{\epsilon_n}(\tilde{v}^j_2),\varphi(x+y^j_{n,2})\rangle+o_n(1)\|\varphi\|
=o_n(1)\|\varphi\|.\endaligned$$
Then
\begin{equation}\label{5.1.4}I'_{\epsilon_n}(w^j_{n,2})\rightarrow0\ \ \text{in}\ \ (H^{\alpha}(\mathbb{R}^N))^*.\end{equation}
Similarly,
$$\tilde{\tilde{I}}_{\epsilon_n}(v^j_{n,1})-\tilde{\tilde{I}}_{\epsilon_n}(w^j_{n,2}(x+y^j_{n,2}))-I_{\nu_2}(\tilde{v}^j_2)\rightarrow0.$$
Then
\begin{equation}\label{5.1.5}I_{\epsilon_n}(w^j_{n,1})-I_{\epsilon_n}(w^j_{n,2})-I_{\nu_2}(\tilde{v}^j_2)\rightarrow0.\end{equation}
 By (\ref{5.1.2}) and (\ref{5.1.5}) we obtain
$$I_{\epsilon_n}(u^j_n)=I_{\epsilon_n}(w^j_{n,2})+I_{V_0}(v^j)+I_{\nu_2}(\tilde{v}^j_2)+o_n(1).$$
Using (\ref{5.1.4}) and (\ref{4.0.0}) we get
$I_{\epsilon_n}(w^j_{n,2})\geq o_n(1).$
Then $I_{\epsilon_n}(u^j_n)\geq c_{V_0}+c_{\nu_2}+o_n(1), $
where $c_{\nu_2}=c_{V(z_0)}$ or $c_{\nu_2}\geq c_{V_\infty}$ using the definition of $\nu_2$.
However,
$I_{\epsilon_n}(u^j_n)=\alpha^j_{\epsilon_n}=c_{V_0}+o_n(1).$
This is a contradiction. Hence, $v^j_n\rightarrow v^j$ in $H^\alpha(\mathbb{R}^N)$.\ \

\textbf{Claim 4.} $y^j_0=a^j.$

Observe that $v^j_n(x)=u^j_n(x+y^j_n)$ and $v^j_n\rightarrow v^j$ in $H^\alpha(\mathbb{R}^N)$. Then for any $\delta>0$ there exists $R_\delta>0$ such that
\begin{equation}\label{7.6.1}\frac{\int_{|x|\geq R_\delta}|u^j_n(x+y^j_n)|^pdx}{\int_{\mathbb{R}^N}|u^j_n(x+y^j_n)|^pdx}\leq\delta.\end{equation}
We firstly show that $y^j_0\in C^j_{l+\nu}$, where $\nu$ is given in (\ref{7.7}).  On the contrary, we have $y^j_0\in \mathbb{R}^N\backslash {C^j_{l+\nu}}$. Since $\epsilon_ny^j_n\rightarrow y^j_0$, we may suppose that for some $i$, $\epsilon_n{(y^j_n)}_i\geq a^j_i+l+\frac{\nu}{2}$ for large $n$.
Then from (\ref{4.4.2}) and (\ref{7.6.1}) we get
$$\aligned h^i_{\epsilon_n}(u^j_n)&=\frac{\int_{|x|\geq R_\delta}\Psi_{\epsilon_n}((x+y^j_n)_i)|u^j_n(x+y^j_n)|^pdx+\int_{|x|\leq R_\delta}\Psi_{\epsilon_n}((x+y^j_n)_i)|u^j_n(x+y^j_n)|^pdx}{\int_{\mathbb{R}^N}|u^j_n(x+y^j_n)|^pdx}\\
&\geq-\frac{2L}{\epsilon_n}\cdot\delta+\Bigl(\frac{a^j_i+l+\frac{\nu}{2}}{\epsilon_n}-R_\delta\Bigr)(1-\delta)
>\frac{a^j_i+l}{\epsilon_n},\endaligned$$
for sufficiently small $\delta>0$, $\nu>\delta$ and $n$ large enough. This is impossible since $\epsilon_nh_{\epsilon_n}(u^j_n)\in C^j_{l}$. Then $y^j_0\in C^j_{l+\nu}$.
Taking into account the condition (V$_2$), (\ref{7.7}) and $V(y^j_0)=V_0$, we deduce
that $y^j_0=a^j$. The proof is complete.  \ \ \ \ $\Box$

Set ${v}^j_{n}:=u^j_n(x+y^j_n)$, where $u^j_n$ is the positive solutions obtained in Lemma \ref{l5.7} and $y^j_n$ is given in (\ref{5.4}). As Lemma \ref{l5.9} and \ref{l5.9.1}, we have the following result.
\begin{lemma}\label{l5.9.2} (i) There exist $C'_j,\delta'_j>0$ such that $C'_j\geq|{v}^j_n|_\infty\geq\delta'_j$ for all $n\in\mathbb{N}$, and $$\lim_{|x|\rightarrow\infty}{v}^j_n(x)=0\ \ \text{uniformly in }n\in\mathbb{N};$$

\noindent(ii) There exists $C_j>0$ such that
$${v}^j_{n}(x)\leq C_j (1+|x|^{N+2\alpha})^{-1},\ \ \forall x\in\mathbb{R}^N,$$
uniformly in $n\in\mathbb{N}$.
\end{lemma}

{\bf Proof of Theorem 1.2.} For any $j\in\{1,...,k\}$, going back
to (\ref{1.1}) with the variable
substitution: $x\mapsto\frac{x}{\epsilon}$, Lemma \ref{l5.7} implies that
 (\ref{1.1}) has $k$ positive solutions for small $\epsilon$.
 Set $\epsilon_n\rightarrow0$. Then arguing as in the proof of Theorem 1.1, from Lemmas \ref{l5.8} and \ref{l5.9.2} we can obtain the desired results. This ends the proof.\ \ \ \ $\Box$

\section{Proof of Theorem 1.3}
In this section, we are going to prove Theorem 1.3.  Recall that
$$\Sigma_\epsilon=\{u\in M_\epsilon: I_\epsilon(u)\leq c_{V_0}+\omega(\epsilon)\},$$
where $\omega(\epsilon)\rightarrow0$ as $\epsilon\rightarrow0$, which is given in Section 4.

\begin{lemma}\label{l6.1}There exists $\bar{\epsilon}>0$ such that for any $\epsilon\in (0,\bar{\epsilon})$,
\begin{equation}\label{6.1}\sup_{u\in \Sigma_\epsilon}\inf_{x\in \cup^k_{j=1}C^j_{{l}/{2\epsilon}}}|h_\epsilon(u)-x|<\frac{l}{2\epsilon}.\end{equation}
\end{lemma}
{\bf Proof}: On the contrary suppose that there exists $\epsilon_n\rightarrow0$ as $n\rightarrow\infty$ such that
\begin{equation*}\sup_{u\in \Sigma_{\epsilon_n}}\inf_{x\in \cup^k_{j=1}C^j_{{l}/{2\epsilon_n}}}|h_{\epsilon_n}(u)-x|\geq\frac{l}{2\epsilon_n}.
\end{equation*}
For each $n$, there exists $u_n\in\Sigma_{\epsilon_n} $ such that
$$\inf_{x\in \cup^k_{j=1}C^j_{{l}/{2\epsilon_n}}}|h_{\epsilon_n}(u_n)-x|=\sup_{u\in \Sigma_{\epsilon_n}}\inf_{x\in \cup^k_{j=1}C^j_{{l}/{2\epsilon_n}}}|h_{\epsilon_n}(u)-x|+o_n(1).$$
In view of $u_n\in \Sigma_{\epsilon_n}\subset M_{\epsilon_n}$, we have
$I_{\epsilon_n}(u_n)\rightarrow  c_{V_0}$ as $n\rightarrow\infty$. By Lemma \ref{l4.5}, there exists $y_n\in\mathbb{R}^N$ such that $v_n(x)=u_n(x+y_n)\rightarrow v_0\neq0$ in $H^\alpha(\mathbb{R}^N)$, and $\epsilon_ny_n\rightarrow y_0\in \mathcal{V}$. By (V$_3$) we know $y_0=a^j$ for some $j=1,...,k$.
Then for large $n$, we have $y_n\in \cup^k_{j=1}C^j_{l/2\epsilon_n}$, where $C^j_{l/2\epsilon_n}$ is given in the first paragraph in Section 5.1. Moreover, since $v_n\rightarrow v_0\neq0$ in $H^\alpha(\mathbb{R}^N)$, for any $\delta>0$, there exists $R_\delta>0$ such that
$\frac{\int_{|x|\geq R_\delta}|v_n|^pdx}{|v_n|^p_p}\leq\delta.$
Then
$$\aligned &|h^i_{\epsilon_n}(u_n)-y^i_n|\\
\leq&\frac{\int_{|x_i|\geq R_\delta}(|\Psi_{\epsilon_n}((x+y_n)_i)|+|y^i_n|)|v_n|^pdx+\int_{|x_i|\leq R_\delta}|\Psi_{\epsilon_n}((x+y_n)_i)-y^i_n||v_n|^pdx}
{|v_n|^p_p}\\
\leq&\delta\cdot \frac{3L}{\epsilon_n}+\frac{\int_{|x_i|\leq R_\delta}|x_i||v_n|^pdx}{|v_n|^p_p}
\leq\delta\cdot \frac{3L}{\epsilon_n}+R_\delta,\quad i=1,...,N,\endaligned$$
for large $n$. Choose $\delta=\frac{l}{12\sqrt{N}L}$, then we have $$\frac{l}{2\epsilon_n}\leq|h_{\epsilon_n}(u_n)-y_n|<\frac{l}{4\epsilon_n}+\sqrt{N}R_{{l}/{(12\sqrt{N}L)}},$$ which is impossible for large $n$. This ends the proof.\ \ \ \ $\Box$

{\bf Proof of Theorem 1.3.} From Lemma \ref{l6.1} it follows that for any $\epsilon\in (0,\bar{\epsilon})$, there holds
$$h_\epsilon(u)\in \cup^k_{j=1}C^j_{{l}/{\epsilon}}, \quad \text{for any } u\in \Sigma_\epsilon.$$ Then $\Sigma_\epsilon\subset \cup^k_{j=1}M^j_\epsilon$ and we can find at least a ground state $u_\epsilon:=w^j_\epsilon\in M^j_\epsilon$ of  (\ref{1.1}) for some $j=1,...,k$, where $w^j_\epsilon$ are the solutions obtained in Theorem 1.2.\ \ \ \ \ $\Box$

\end{document}